\documentclass{amsart}
\usepackage[letterpaper, portrait, margin=2.5cm]{geometry}

\usepackage{amsmath,amssymb,amsthm}
\usepackage{mathtools}
\usepackage{xcolor,colortbl}
\usepackage[all,cmtip]{xy}
\input diagxy

\usepackage{hyperref}
\usepackage{tikz-cd}
\usepackage{caption}
\usepackage{enumitem}

\newcommand{\Ab}{\mathbb A}

\newcommand{\Gb}{\mathbb G}

\newcommand{\Nb}{\mathbb N}

\newcommand{\Pb}{\mathbb P}
\newcommand{\Qb}{\mathbb Q}

\newcommand{\Vb}{\mathbb V}

\newcommand{\Zb}{\mathbb Z}

\newcommand{\bd}{\mathbf{d}}
\renewcommand{\bf}{\mathbf{f}}
\newcommand{\bg}{\mathbf{g}}
\newcommand{\bm}{\mathbf{m}}

\newcommand{\bh}{\mathbf{h}}

\newcommand{\mcc}{\mathcal{C}}

\newcommand{\mcf}{\mathcal{F}}
\newcommand{\mcg}{\mathcal{G}}

\newcommand{\mco}{\mathcal{O}}
\newcommand{\mcp}{\mathcal{P}}

\newcommand{\mcr}{\mathcal{R}}

\newcommand{\mct}{\mathcal{T}}
\newcommand{\mcu}{\mathcal{U}}

\newcommand{\mcx}{\mathcal{X}}

\newcommand{\mcz}{\mathcal{Z}}

\newcommand{\lbm}{\left[ \begin{matrix}}
\newcommand{\rem}{\end{matrix} \right]}

\DeclareMathOperator{\GL}{GL}
\DeclareMathOperator{\PGL}{PGL}

\DeclareMathOperator{\grass}{Gr}

\newcommand{\floor}[1]{\left\lfloor {#1} \right\rfloor}
\newcommand{\Mt}{M_{\operatorname{type}}}
\newcommand{\Md}{M_{\operatorname{deg}}}

\DeclareMathOperator{\Perm}{Perm}

\newcommand{\into}{\hookrightarrow}

\DeclareMathOperator{\spec}{Spec}

\DeclareMathOperator{\red}{red}

\DeclareMathOperator{\characteristic}{char}
\newcommand{\Sol}{\mathrm{Sol}}

\DeclareMathOperator{\Sym}{Sym}

\DeclareMathOperator{\codim}{codim}

\newcommand{\density}[2]{f_{#1}(#2)}
\newcommand{\existence}[2]{g_{#1}(#2)}
\newcommand{\param}[2]{\Ab_{#1}(#2)}
\newcommand{\polar}[2]{{#1}^{[#2]}}

\theoremstyle{definition}
\newtheorem{definition}{Definition}[section]
\newtheorem{problem}[definition]{Problem}
\newtheorem{question}[definition]{Question}
\newtheorem{example}[definition]{Example}

\newtheorem{convention}[definition]{Convention}

\theoremstyle{remark}
\newtheorem{remark}[definition]{Remark}

\theoremstyle{plain}
\newtheorem{prop}[definition]{Proposition}
\newtheorem{lemma}[definition]{Lemma}
\newtheorem{theorem}[definition]{Theorem}
\newtheorem{corollary}[definition]{Corollary}

\title{Solvable points on intersections of quadrics, cubics, and quartics}
\author{Claudio G\'{o}mez-Gonz\'{a}les}
\email{cgonzales@carleton.edu}
\address{Department of Mathematics \& Statistics, Carleton College}

\author{Jesse Wolfson}
\email{wolfson@uci.edu}
\address{Department of Mathematics, University of California-Irvine}

\thanks{The first author was supported in part by NSF Grant DMS-2418943. The second author was supported in part by NSF Grants DMS-1944862 and DMS-2506184.}

\begin{document}

\begin{abstract} Let $k$ be a field of characteristic not 2 or 3.  We establish polynomial bounds on the ambient dimension $N$, such that if $N$ is at least as big as this bound, every intersection $X\subset\Pb^N$ of quadrics, cubics and quartics has a solvable point, and a general one has a dense collection of solvable points, i.e. points in $X(k^{\Sol})$ where $k^{\Sol}/k$ is a solvable closure. Our method connects the classical theory of polar hypersurfaces to Fano varieties $\mcf(j,X)$ of linear subspaces on $X$, and we use this to obtain improved control on the arithmetic of $\mcf(j,X)$.
\end{abstract}

\maketitle

\section{Introduction}
Let $k$ be a field and let $X$ be an irreducible variety over $k$. The study of rational points $X(k)$ is deep and well-established. By contrast, much less has been written about solvable points on varieties, i.e., points $X(k^{\Sol})$, where $k^{\Sol}/k$ is a solvable closure. For example, if $k=\Qb$, consider the following.
\begin{question}\label{q:innocuous}
    Does every geometrically irreducible $\Qb$-variety $X$ admit a dense collection of solvable points?
\end{question} 
Note that there are no local obstructions to solvable points, because the absolute Galois group of $\Qb_p$ is pro-solvable. Moreover, because $\Qb^{\Sol}$ has cohomological dimension $1$, the known global cohomological obstructions, e.g., Brauer--Manin, vanish, and we find ourselves in a situation where the answer to Question~\ref{q:innocuous} could be, for all we know, ``Yes.''  

Since the turn of this century, a small number of papers, beginning with work of P\'al \cite{Pal2004}, have been written on this subject, showing that varieties whose geometry is sufficiently simple admit solvable points. The deepest of these results are due to \c{C}iperiani-Wiles \cite{CiperianiWiles2008}, who use modularity and systems of Heegner points to show that smooth genus $1$ curves over $\Qb$ satisfying mild assumptions have a solvable point. Others have considered the analogous question over arbitrary fields of characteristic 0; e.g., Wooley \cite{Wooley2016} applies a ``diagonalization'' method introduced by Brauer \cite{Brauer1945} to show that when the degree of a smooth, geometrically integral projective variety is small relative to its dimension, then it admits a solvable point. Unfortunately, the methods do not address density and the bounds on the dimension that Wooley achieves are enormous.

In this paper, we consider a method due to Sylvester \cite{Sylvester1887} for obtaining solutions to systems of homogeneous equations in terms of polynomials of lower degrees if the number of variables in the system is sufficiently large. Sylvester's ``method of obliteration'' is remarkably general, yet appears to have been largely forgotten in modern sources: see \cite{Heberle2021} for a hands-on introduction, and \cite{HeberleSutherland2023} for a treatment in terms of polar cones that we generalize here. In this paper, we formalize a variation of Sylvester's obliteration algorithm to address the existence of solvable points on varieties. Indeed, for intersections of hypersurfaces of low degree, we show that varieties of sufficiently large dimension (though significantly smaller than previously described in the literature) have many solvable points.

\begin{theorem}[Existence and density of solvable points]\label{theorem:solvable_points}
Let $k$ be a field of characteristic not $2$ or $3$, and let $m_2, m_3, m_4 \in \Nb$. There exists a polynomial $p \in \Qb[m_2,m_3,m_4]$ independent of $k$ with non-negative coefficients and total degree $8$ such that, whenever $N \geq p(m_2,m_3,m_4)$, we have:
\begin{enumerate}[label=(\alph*)]
    \item\label{i:existence_pts} Any intersection $X\subseteq\Pb^N$ of $m_2$ quadrics, $m_3$ cubics and $m_4$ quartics, has $X(k^{\Sol}) \not = \emptyset$.
    \item\label{i:density_pts} For \emph{general} intersections $X\subseteq\Pb^N$ of $m_2$ quadrics, $m_3$ cubics and $m_4$ quartics, as in Definition \ref{definition:general_property}, $X(k^{\Sol}) \subseteq X$ is dense.
\end{enumerate}
\end{theorem}

\noindent ``General'' has its usual meaning here: there is a non-empty Zariski open subset $U$ of the affine space of $(m_2+m_3+m_4)$-tuples of forms of the given degrees, such that the conclusion holds for every tuple in $U$ with coefficients in $k$. In contrast, part~\ref{i:existence_pts} holds for \emph{every} intersection of the given type.

\begin{remark}\label{remark:counterexample}
Part~\ref{i:density_pts} of Theorem~\ref{theorem:solvable_points} cannot be strengthened to all intersections $X$, as was first pointed out to us by Boris Alexeev using a pre-release version of Astra, and then later by David Savitt using GPT-5.6 Sol (Ultra High), as well as by the anonymous referee. Indeed, suppose $X\subseteq\Pb^n$ has $X(k^{\Sol})=\emptyset$ (e.g., the zero-dimensional scheme $\Vb(yz-x^2,\ y^3-xz^2+z^3)\subseteq\Pb^2_\Qb$, whose points are defined over $\Qb(\alpha)$ with $\alpha^6-\alpha+1=0$). Then, for any $N>n$, the same equations define a cone $Y\subseteq\Pb^N$ with vertex a linear subspace $\Pb^{N-n-1}$ (an intersection of the same type) and projection from the vertex maps the rest of $Y$ onto $X$. Hence a solvable point of $Y$ outside the vertex would give a solvable point of $X$, so $Y(k^{\Sol})$ is contained in the vertex and not dense in $Y$ for any $N$. 
\end{remark}


Moreover, our method of proof makes apparent a more general result regarding $j$-dimensional planes on varieties. We write $\mcf(j,X)$ for the scheme of $j$-planes on $X$ and $q \in \Qb[j,m_2,m_3,m_4]$ for a total degree $8$ polynomial with non-negative coefficients given explicitly in Appendix~\ref{section:select_bounds}.

\begin{theorem}[Universal existence and general density of solvable $j$-planes]\label{theorem:solvable_j-planes}
For any $j \in \Nb$ and with $k$ and $m_2, m_3, m_4 \in \Nb$ as before, whenever $N \geq q(j,m_2,m_3,m_4)$,
\begin{enumerate}[label=(\alph*)]
    \item\label{i:existence_planes} Any intersection $X\subseteq\Pb^N$ of $m_2$ quadrics, $m_3$ cubics and $m_4$ quartics, has $\mcf(j,X)(k^{\Sol}) \not = \emptyset$.
    \item\label{i:density_planes} For general intersections $X\subseteq\Pb^N$ of $m_2$ quadrics, $m_3$ cubics and $m_4$ quartics, $\mcf(j,X)(k^{\Sol}) \subseteq \mcf(j,X)$ is dense.
\end{enumerate}
\end{theorem}

While $p$ and $q$ are ungainly to typeset in full generality in an introduction (see Appendix~\ref{section:select_bounds} for the explicit formula), specific cases make the point clear:

\begin{corollary}\label{corollary:points_on_pure_intersections} For a general intersection $X \subseteq \Pb^N$ of $m \geq 1$ hypersurfaces, $X(k^\Sol) \subseteq X$ is dense if:
\begin{enumerate}[label=(\alph*)]
    \item each hypersurface is of degree $2$ and $N \geq \tfrac{1}{2} (m+1)^2$.
    \item each hypersurface is of degree $3$ and $N \geq \tfrac{1}{2^3} (m+1)^4$.
    \item each hypersurface is of degree $4$ and $N \geq \tfrac{1}{2^7} (m+1)^8$.
\end{enumerate}
\end{corollary}
\noindent These estimates eschew sharpness (with regards to the methods described in this paper) for tidiness. We include them to compare with \cite[Theorem 2]{Wooley1998} which guarantees a solvable point for systems of $m$ forms of degree $d$ in $\Pb^N$ for large enough $N$; in particular, when $N \geq 2m^2$, $N \geq 4m^4$, and $N \geq 16m^8$ for the respective cases of quadrics, cubics, and quartics. We emphasize that these bounds only guarantee the existence of \emph{a} solvable point, not the density of such, as emphasized in \cite[Problem (b)]{Wooley1999}.\footnote{As Wooley informed us, using \cite[Lemma 2.1]{Wooley1998} and $\varphi=1$, one can obtain the same constants as in Corollary~\ref{corollary:points_on_pure_intersections} for existence (though not density).} 

Indeed, even the bounding polynomial $q$ is not the sharpest bound obtainable through these methods---the formulas in Section~\ref{subsection:iterative_formulas} involve floor functions which admit additional savings. See Appendix \ref{section:select_bounds} for select tables of our tightest bounds on $N$ as in Theorem \ref{theorem:solvable_j-planes}.

\medskip
Sylvester's obliteration algorithm produces a point of $X=Y\cap Z$ by first producing a $j$-plane $\Lambda\subseteq Z$ and then solving the low-degree equations of $Y$ on $\Lambda\cong\Pb^j$; iterating, one trades the equations of $X$ for equations of lower degree on a smaller projective space, at the cost of ambient dimension. The $j$-planes themselves are produced inductively: a $(j-1)$-plane $\Lambda\subseteq X$ and a point of the \emph{polar model} $\mcf_\Lambda(j,X)$---a subscheme of $\Pb^{N-j}$ cut out by the polars of the equations of $X$ along $\Lambda$ (Theorem~\ref{thm:type_of_fano})---together span a $j$-plane in $X$. With existence of a solvable point as the inductive hypothesis, this argument applies to every intersection and yields part~\ref{i:existence_planes} of Theorem~\ref{theorem:solvable_j-planes} (Section~\ref{subsection:existence}).

For density, each step of the induction passes a dense set of solvable points from the base of a correspondence to its total space. This is only possible if every irreducible component of the total space dominates the base (Lemma~\ref{lemma:fibration}). The cones of Remark~\ref{remark:counterexample} show how this may fail: 
the planes through the vertex which meet $Y$ away from the vertex form a proper closed subset of $\mcf(j,Z)$ along which the fibers $\Lambda\cap Y$ have positive dimension, and the component of the total space lying over this locus---the one that dominates $Y\cap Z$---does not dominate $\mcf(j,Z)$.
We therefore work throughout with \emph{general} intersections, in the explicit sense of a dense Zariski open subset of the parameter space $\param{N}{\bd}$ of defining equations. Two facts make the induction possible. First, for general $X$ the dominance condition can be verified, either because $\mcf(j,X)$ is irreducible by a theorem of Debarre and Manivel \cite{DebarreManivel1998}, or by a dimension count on the locus of planes meeting $Y$ excessively (Lemma~\ref{lemma:dominance}). Second, the equations of the polar model $\mcf_\Lambda(j,X)$ depend explicitly on the equations of $X$ containing $\Lambda$ (Theorem~\ref{thm:linear_model}). In particular, this implies that a suitably general $X$ has suitably general polar models along a dense set of its solvable $(j-1)$-planes, and the inductive hypothesis can be applied to them.

\bigskip
Our approach to Sylvester's obliteration algorithm uses a classical perspective on varieties $\mcf(j,X)$ of linear subspaces on intersections of hypersurfaces which appears to be absent from contemporary treatments such as \cite{DebarreManivel1998}. Building off the theory of iterated polar cones from Sutherland's thesis \cite{Sutherland2021}, we introduce the ``total $j$-polar space'' $\mcc^j(X)$ as an incidence variety living over $\mcf(j,X)$. In this context, Sutherland's work allows us to get finer control of the arithmetic of these varieties than the typical approach via Pl\"ucker coordinates and the geometry of the Grassmannian. We refer the reader to Section~\ref{s:geom} for this treatment, which we believe may be of interest beyond the application to solvable points.

We believe that the study of solvable points is a natural question from a number of perspectives. In addition to the relation to rational points noted at the beginning, the existence of dense collections of solvable points is closely related to classical questions such as Hilbert's 13th Problem on how to find the simplest formula for a polynomial of one variable. As another example, an elementary argument \cite[Proposition 11.1.1]{FriedJarden2008} shows that Question~\ref{q:innocuous} is equivalent to the following open problem.

\begin{problem}\cite[Problem 11.5.9(a)]{FriedJarden2008}
    Is $\Qb^{\Sol}$ a pseudo-algebraically closed (PAC) field? 
\end{problem}

\noindent See \cite[Remark 11.5.10]{FriedJarden2008} for a discussion of this problem in connection with classical questions on Riemann surfaces. 

Given the above, we find the sparsity of the contemporary literature on solvable points surprising. Indeed, the following forms a complete bibliography to our knowledge of engagements with solvable points in the present century \cite{Pal2004,CiperianiWiles2008,Pal2013,Wooley2016,Rawson2024}. We hope that this paper and possible sequels might contribute to more mathematicians engaging with this topic.

\subsection*{Statement on AI Use}
The original version of this paper appeared on the arXiv in the summer of 2025, and was written without AI. As discussed in Remark~\ref{remark:counterexample}, in mid-August 2026, Boris Alexeev, using a pre-release version of Astra, found an error in our original statement of Theorem~\ref{theorem:solvable_points}, which was independently found by David Savitt a few weeks later using GPT-5.6 Sol. After we identified the present repair (i.e. replacing the original ``for all'' with the present ``for general''), we used Claude Science (Fable 5.1) to assist in the revision process, giving it our proposed repair, pointing it to the necessary places in the literature, and asking it to track the logical dependencies of the fix and to draft a detailed implementation, which we then substantially reworked to produce the revised text.  We used Claude Science (Fable and Opus) for a final review before resubmitting. The authors themselves generated the arguments, checked each detail, and take sole responsibility for the paper's correctness.

\subsection{Acknowledgements}
The authors thank Alexander Sutherland for introducing them to the perspective of polar cones and their applications.  We thank Reginald Anderson, Jordan Ellenberg, James Rawson, Zinovy Reichstein, and Trevor Wooley for generous comments on an early draft. We thank Boris Alexeev for alerting us to an error in the original version of Theorem~\ref{theorem:solvable_points}; David Savitt independently found the error and alerted us a few weeks later, as did the anonymous referee. We thank the anonymous referee for additional substantive comments that improved the paper. 

\section{Preliminaries}\label{s:prelim}

\subsection{Notation}

We recall the necessary geometric background, specify notations, and give a brief treatment of polar cones.

\begin{definition}\label{definition:type}
For any $n$, given hypersurfaces $X_{1,1},\dots,X_{1,m_1},X_{2,1},\dots,X_{2,m_2},\dots,X_{n,1},\dots,X_{n,m_n} \subset \Pb^N$, where each $X_{i,j}$ is of degree $i$, we say that their (scheme theoretic) intersection $X:=\bigcap_{i,j} X_{i,j}$ is of \emph{type} $\bm = (m_1,\dots,m_n) \in \Nb^n$. More generally, for $X$ a (scheme theoretic) intersection of hypersurfaces of type $\bm$, we say that $X^{\red}$ is also of type $\bm$.
\end{definition}

Our interest in this paper is in working over non-algebraically closed fields $k$.
\begin{convention}\label{convention:surjective}
    By a {\em surjective} map of $k$-schemes $f: X\to Y$, we mean one that is surjective on geometric points, i.e. $f: X(\Omega)\to Y(\Omega)$ is a surjective map of sets, for any algebraically closed extension $\Omega/k$.
\end{convention}

The main Theorems of this paper concern intersections of type $(0,m_2,m_3,m_4)$. Note that our notation is different from \cite{Sutherland2021}, where such a type would be written as $\big( \begin{smallmatrix}
    4 & 3 & 2 & 1 \\
    m_4 & m_3 & m_2 & 0
\end{smallmatrix} \big)$ or, more generally, as \[ \begin{pmatrix}
    n & \cdots & 2 & 1 \\
    m_n & \cdots & m_2 & m_1
\end{pmatrix}. \]
We will also make use of the notation of Debarre--Manivel \cite{DebarreManivel1998}, which would encode the aforementioned intersection by a tuple
\begin{equation}\label{eqn:degree_from_type}
    \bd = (\underbrace{1,\dots,1}_{m_1 \text{ many}},\underbrace{2,\dots,2}_{m_2 \text{ many}},\dots\dots,\underbrace{n,\dots,n}_{m_n \text{ many}}),
\end{equation}
referring to the intersection $X$ as a \emph{variety defined by equations of degree $\bd$}. Note that such a tuple need not be written in ascending order, though we have done so here for clarity. This notation has the advantage of coinciding with a tuple of polynomials $\bf$ cutting out $X$ when the specific $f_{i,j}$ cutting out $X_{i,j}$ are known. For the remainder of the paper, we use $\bm$ for type, $\bd$ for degrees, and $\bf$ for the associated polynomials.

\medskip

\begin{remark}\label{remark:type}
If $X\subseteq \Pb^N$ is an intersection of hypersurfaces, neither its type, nor the type of $X^{\red}$, is unique. In general, we are only careful to speak of \emph{the} type when referring to a specific (perhaps unnamed) set of polynomials cutting out $X$.
\end{remark}

\begin{remark}\label{remark:reduced}
    For any scheme $X$, and any field $L$, $X(L)=X^{\red}(L)$. As our primary interest in this paper is in rational points, our underlying interest is thus in $X^{\red}$, even if we distinguish between $X$ and $X^{\red}$ for the sake of geometry.
\end{remark}

\subsection{Combinatorics}

In this section we introduce two monoids $\Mt$ and $\Md$, whose underlying sets are both the collection of all compactly supported sequences $C_c(\Nb_{\geq 1}, \Nb) = \varinjlim_{n \geq 1} \Nb^n$. Their binary operations, and other iterative operations we will define here, are understood as operations on types and degrees as the result of procedures carried out on underlying polynomials. We will often write elements of $\Mt$ and $\Md$ with rightmost zeroes suppressed, i.e., simply as $(m_1,\dots,m_n)$ and $(d_1,\dots,d_r)$, respectively, with the understanding that $m_n > 0$ and $d_r > 0$. We also verify purely arithmetic results regarding these operations in service of cleaner exposition for the geometric arguments to come.

\begin{definition} 
We define $\Mt$ as the monoid $\left( C_c(\Nb_{\geq 1}, \Nb), + \right)$, where $+$ is entrywise addition.
 For non-zero sequences $\bd = (d_1,\dots,d_r)$ and $\bd' = (d_1',\dots,d_s')$, where $d_r \not = 0 \not = d_s'$, we define
\[
	\bd \oplus \bd' = (d_1,\dots,d_r, d_1',\dots,d_s'),
\]
and accordingly say $(0,\dots) \oplus \bd = \bd = \bd \oplus (0,\dots)$. We set $\Md = \left( C_c(\Nb_{\geq 1}, \Nb), \oplus \right)$.
\end{definition}

Note that there is a natural homomorphism $t: \Md \to \Mt$ sending a tuple of degrees to the underlying type as in Definition \ref{definition:type}. Moreover, given a type $\bm$, there is a unique weakly increasing degree tuple $\bd$ such that $t(\bd) = \bm$, namely \eqref{eqn:degree_from_type}; we refer to it as \emph{the} degree tuple of $\bm$.

\begin{example}\label{ex:plus}
    Let $X$ and $Y$ be subschemes in $\Pb^N$ cut out by collections of homogeneous polynomials with degree tuples given by $\bd$ and $\bd'$, and types $\bm = t(\bd)$ and $\bm' = t(\bd')$. Then $X\cap Y$ is cut out by equations of degree $\bd\oplus\bd'$ and type $\bm+\bm'$.
\end{example}

\begin{definition}\label{d:type}
Given a type $\bm = (m_1, \dots, m_n) \in \Mt$, we write
\[
    \bm^1 \coloneqq \Bigg(\sum_{i=1}^n m_i, \dots, \underbrace{\sum_{i=\ell}^n m_i}_{\ell\text{th entry}}, \dots, m_n \Bigg).
\]
Moreover, for any $j \in \Nb_{\geq 1}$, we define $\bm^{j+1} \coloneqq (\bm^j)^1$ and $\bm^0 \coloneqq \bm$. For $\bd = (d_1, \dots, d_r) \in \Md$, we set
\[
    \bd^1 \coloneqq \bigoplus_{i=1}^r (1, 2, \dots, d_i)
\]
and again $\bd^{j+1} = (\bd^j)^1$ for all $j \in \Nb_{\geq 1}$, $\bd^0=\bd$.
We also make use of the $\ell^1$ norm
\[
    |\bm|=\sum_{i=1}^\infty m_i.
\]
\end{definition}

\begin{lemma}\label{l:comb}
Let $t: \Md \to \Mt$ be the homomorphism sending a tuple of degrees to the underlying type as in Definition \ref{definition:type} and let $j \in \Nb_{\geq 1}$. Then:
\begin{enumerate}[label=(\alph*)]
\item\label{i:monoid_endo} The assignments $\bm\mapsto \bm^j$ and $\bd\mapsto \bd^j$ are monoid endomorphisms on $\Mt$ and $\Md$ respectively.
\item\label{i:type_commutes} For any $\bd \in \Md$, we have $t(\bd^j) = t(\bd)^j$.
\item\label{i:entries} For any $\bm \in \Mt$, we have \[
    \bm^j = \left( \sum_{i=1}^n \binom{j+i-2}{i-1} m_i, \dots,
    \sum_{i=\ell}^n \binom{j+i-\ell-1}{i-\ell} m_i, \dots, m_n \right).
\]
\item\label{i:norm} For any $\bm\in\Mt$, we have 
\[ |\bm^j| = \sum_{i=1}^n \binom{i+j-1}{j} m_i. \] 
In particular, $\bm^{j-1} \leq \bm^{j}$ entrywise and $|\bm^{j-1}| \leq |\bm^j|$.
\end{enumerate}
\end{lemma}
\begin{proof}
The identities $(\bm+\bm')^j=\bm^j+\bm'^j$ and $(\bd\oplus\bd')^j = \bd^j \oplus \bd'^j$ follow from the $j=1$ cases and the inductive nature of our definitions.

We prove the second claim by induction on $j$. In light of the previously established identities, it suffices to demonstrate the equality in the case of a singleton degree $\bd = (d)$. The straightforward calculation
\[
    t\left( (d)^1 \right) = t(1,2,\dots,d) = (1,1,\dots,1) = (0,\dots,0,1)^1 = t(d)^1
\]
establishes the base case for all $\bd$. For the inductive step, having established the formula for some fixed $j$, we simply have:
\[
    t\left( \bd^{j+1} \right) = t\left((\bd^j)^1\right) = t\left(\bd^j\right)^1 = \left(t(\bd)^j\right)^1 = t(\bd)^{j+1}.
\]
The third claim also follows by induction on $j$, where the base case is a tautology. For the inductive step, we consider the $\ell$th entry of $(\bm^j)^1$:
\[
\sum_{s = \ell}^n \sum_{i=s}^n \binom{j+i-s-1}{i-s} m_i = \sum_{i=\ell}^n \sum_{s' = 0}^{i-\ell} \binom{s'+j-1}{s'} m_i = \sum_{i=\ell}^n \binom{j+i-\ell}{i-\ell} m_i
\]
where the first equality interchanges sums via setting $s' = i-s$ and the second is Chu's identity.

For the fourth claim, summing the $\ell$th entries in \ref{i:entries} over $\ell=1,\dots,i$ for fixed $i$ and applying Chu's identity once more gives 
\[ \sum_{\ell=1}^{i}\binom{j+i-\ell-1}{i-\ell}=\binom{j+i-1}{i-1}=\binom{i+j-1}{j}. \]
Monotonicity in $j$ is clear from \ref{i:entries}, since the binomial coefficients $\binom{j+i-\ell-1}{i-\ell}$ are weakly increasing in $j$.
\end{proof}

\begin{example}\label{example:type_shift}
For $\bm=(0,m_2,m_3,m_4)$ and $j\geq 0$, Lemma~\ref{l:comb} gives
\begin{align*}    
    \bm^j &=\Big(jm_2+\tbinom{j+1}{2}m_3+\tbinom{j+2}{3}m_4, m_2+jm_3+\tbinom{j+1}{2}m_4, m_3+jm_4, m_4\Big), \\
    |\bm^j| &=(j+1)m_2+\tbinom{j+2}{2}m_3+\tbinom{j+3}{3}m_4.
\end{align*}
These two expressions will recur throughout Section~\ref{s:solvable}: the first entry of $\bm^j$ counts linear forms that will be absorbed into the ambient dimension, and $|\bm^j|+j$ is the threshold beyond which every point of an intersection of type $\bm$ lies on a $j$-plane (Corollary~\ref{c:tf3}).
\end{example}

\subsection{Rational points in fibrations}\label{subsection:fibrations}

Our strategy for producing rational points is inductive: rational points on a base, together with rational points on fibers, produce rational points on a total space. Throughout this subsection, $K$ is an arbitrary field, a variety is a separated scheme of finite type over $K$, irreducible components carry their reduced structure, and $\dim\emptyset=-\infty$.

\begin{lemma}\label{lemma:fibration}
Let $\pi: E \to B$ be a morphism of $K$-varieties. Assume that:
\begin{enumerate}[label=(\alph*)]
\item\label{item:dense-base} there is a subset $\Sigma \subseteq B(K)$ that is Zariski dense in $B$; 
\item\label{item:dense-fibers} for every $b \in \Sigma$, the fiber $\pi^{-1}(b)$ has a Zariski dense set of $K$-points; and 
\item\label{item:dominate-components} every irreducible component of $E$ dominates an irreducible component of $B$. 
\end{enumerate}
Then $E(K) \subseteq E$ is Zariski dense.
\end{lemma}
\begin{proof}
Let $E'$ be an irreducible component of $E$, which by \ref{item:dominate-components} dominates some irreducible component $B'$ of $B$. Let $V \subseteq E$ be an arbitrary non-empty open subset intersecting $E'$; we must show $V$ contains a $K$-point.

Shrinking $V$, we may assume that $V \subseteq E'$ and that $V$ is disjoint from the other irreducible components of $E$. In this case, $V$ is dense in $E'$, so $\pi(V)$ is dense in $\pi(E')$ and hence dense in $B'$. Moreover, because $\pi(V)$ is a constructible set (by Chevalley's theorem), it contains a non-empty open subset of $B'$. Shrinking this subset to avoid the components of $B$ not equal to $B'$, we obtain a nonempty open $W \subseteq B$ with $W \subseteq \pi(V)$.

Because $\Sigma$ is dense in $B$ by \ref{item:dense-base}, there is some $b \in \Sigma \cap W$. Since $b \in \pi(V)$, the set $V \cap \pi^{-1}(b)$ is a non-empty open subset of $\pi^{-1}(b)$. Finally, by \ref{item:dense-fibers}, the $K$-points of this fiber are dense, so $V \cap \pi^{-1}(b)$ contains a point of $\pi^{-1}(b)(K) \subseteq E(K)$.
\end{proof}

\begin{corollary}\label{corollary:density-bundle}
Let $\pi: E \to B$ be a Zariski locally trivial fiber bundle of $K$-varieties whose fiber $\Phi$ is geometrically irreducible with $\Phi(K) \subseteq \Phi$ dense (e.g., $\Phi=\Ab^n$, $\Pb^n$, or a non-empty open subset of either). If we have $B(K) \subseteq B$ dense, then $E(K) \subseteq E$ is dense.
\end{corollary}
\begin{proof}
Hypotheses \ref{item:dense-base} and \ref{item:dense-fibers} of Lemma \ref{lemma:fibration} hold with $\Sigma=B(K)$, and hypothesis \ref{item:dominate-components} holds because the geometric irreducibility of the fibers ensures that the irreducible components of $E$ are the preimages of the irreducible components of $B$.
\end{proof}

Finally, we record a criterion for verifying the dominance hypothesis \ref{item:dominate-components} of Lemma \ref{lemma:fibration} by a dimension count. In our applications, $E$ will be cut out of a projective or affine bundle $P\to B$ by one or two defining equations, and the locus in $B$ over which the fibers of $E\to B$ are larger than expected will have small dimension.

\begin{lemma}\label{lemma:dominance}
Let $\varphi: P \to B$ be a morphism of $K$-varieties, and let $a,c \in \Nb$. Suppose that for every irreducible component $B'$ of $B$, the preimage $\varphi^{-1}(B')$ is irreducible of dimension $\dim B' + a$. Let $E \subseteq P$ be the common zero locus of $c$ sections of line bundles on $P$, and let $\pi \coloneqq \varphi|_E: E \to B$. For any integer $s$, define the constructible subset of $B$:
\[
  B_{\ge s} \coloneqq \{ b \in B \mid \dim \pi^{-1}(b) \ge s \}.
\]
Suppose that for every irreducible component $B'$ of $B$ and every $s > a - c$, we have
\begin{equation}\label{eq:domcount}
  \dim(B_{\ge s} \cap B') + s < \dim B' + a - c.
\end{equation}
Then every irreducible component of $E$ dominates an irreducible component of $B$.
\end{lemma}
\begin{proof}
Let $E'$ be an irreducible component of $E$. Because $P$ is the union of the preimages of the irreducible components of $B$, $E'$ is contained in $\varphi^{-1}(B')$ for some irreducible component $B'$ of $B$.

We first establish a lower bound on the dimension of $E'$. By hypothesis, $E \cap \varphi^{-1}(B')$ is cut out from the irreducible variety $\varphi^{-1}(B')$ by $c$ sections of line bundles. By Krull's theorem, every irreducible component of this intersection has codimension at most $c$ in $\varphi^{-1}(B')$, and so
\begin{equation}\label{eq:krull}
  \dim E' \ge \dim \varphi^{-1}(B') - c = \dim B' + a - c.
\end{equation}

By definition, the fibers of $\pi|_{E'}$ over the constructible set $(B_{\ge s} \setminus B_{\ge s+1}) \cap B'$ have dimension at most $s$. Therefore, for every $s>a-c$,
\[
  \dim \big(E'\cap \pi^{-1}((B_{\ge s} \setminus B_{\ge s+1}) \cap B')\big) \leq \dim(B_{\ge s} \cap B') + s < \dim B' + a - c,
\]
where the strict inequality is our hypothesis \eqref{eq:domcount}. Comparing with \eqref{eq:krull}, these preimages have strictly smaller dimension than $E'$. Since fiber dimensions are bounded, $E'\cap\pi^{-1}(B_{\ge a - c + 1})$ is a finite union of such sets, and thus cannot be all of $E'$.
Therefore $V \coloneqq E' \setminus \overline{\pi^{-1}(B_{\ge a - c + 1})}$ is a non-empty open subset of $E'$, and every fiber of $\pi|_V$ has dimension at most $a - c$.

Now, because $V$ is dense in $E'$, we have $\overline{\pi(V)} = \overline{\pi(E')}$. Since the fibers of $\pi|_V$ have dimension at most $a - c$, we obtain the upper bound
\[
  \dim E' = \dim V \le \dim \overline{\pi(E')} + a - c.
\]
Combining with \eqref{eq:krull} yields $\dim \overline{\pi(E')} \ge \dim B'$ and, by the irreducibility of $B'$, we must have $\overline{\pi(E')} = B'$.
\end{proof}

\begin{remark}\label{remark:dominance_krull}
The proof only uses the consequence \eqref{eq:krull} of the hypothesis on $E$; Lemma \ref{lemma:dominance} therefore holds for any closed subset $E \subseteq P$ such that every irreducible component of $E \cap \varphi^{-1}(B')$ has dimension at least $\dim B' + a - c$. Note also that if every fiber of $\pi$ has dimension exactly $a-c$, then $B_{\geq s}=\emptyset$ for $s>a-c$ and the hypothesis \eqref{eq:domcount} is vacuous.
\end{remark}

\section{Geometric constructions}\label{s:geom}

We write $\grass(j,N)$ for the Grassmannian of projective $j$-planes in $\Pb^N$; recall the correspondence with the usual Grassmannians ($\mathrm{G}$, of subspaces of vector spaces): $\grass(j,N) \cong \mathrm{G}(j+1,N+1)$. We denote the tautological (sub-)bundle by
\[
    S_j\to \grass(j,N).
\]
In addition, given a closed subscheme $X \subseteq \Pb^N$, we denote its scheme of $j$-planes by
\[ \mcf(j,X) \coloneqq \{ \Lambda \in \grass(j,N): \Lambda \subseteq X \}. \]
If $x \in X$, we also write $\mcf_x(j,X)$ for the $j$-planes in $X$ containing $x$; more generally, given an $\ell$-plane $\Lambda\subset X$, we write $\mcf_\Lambda(j,X)\subset \mcf(j,X)$ for the space of $j$-planes in $X$ which contain $\Lambda$, and $\grass_{\Lambda}(j,N)=\mcf_\Lambda(j,\Pb^N)$. Note that the restriction
\[
    S_j|_{\grass_\Lambda(j,N)}\to\grass_\Lambda(j,N)
\]
admits a trivial sub-bundle $\grass_\Lambda(j,N)\times \Lambda\subset S_j|_{\grass_\Lambda(j,N)}$, and thus likewise for $S_j|_{\mcf_\Lambda(j,X)}$. By abuse of notation, we denote this trivial sub-bundle by
\[
    \Lambda\subset S_j|_{\grass_\Lambda(j,N)}.
\]
Given $x_0,\ldots,x_j\in \Pb^N$, we denote the linear subspace of $\Pb^N$ they span by $\Lambda(x_0,\ldots,x_j)$; more generally, for a linear subspace $\Lambda$ and a point $y$, we write $\langle \Lambda, y\rangle$ for their span. Lastly, we denote the subscheme defined by $f_1,\dots,f_\ell \in k[z_0,\dots,z_N]$ using the traditional $\Vb(f_1,\dots,f_\ell) \subseteq \Pb^N$.

\subsection{Incidence varieties of linear subspaces}

Let $\bd=(d_1,\ldots,d_r)\in\Md$ be a tuple of positive integers. For the remainder of this section, let $X\subset \Pb^N$ be a subscheme cut out by a collection of homogeneous polynomials $\bf$ of degree $\bd$. 
For $j\ge 0$, recall (e.g., \cite[Remarques 2.10]{DebarreManivel1998}) that $\mcf(j,X)\subset \grass(j,N)$ can be defined as the vanishing locus of a section of a vector bundle. More precisely, letting $S_j\to \grass(j,N)$ be the tautological bundle as above, the tuple of polynomials $\bf$ defines a section of
\begin{equation}\label{eqn:fano_section}
    \Sym^{\bd} S_j^\ast \coloneqq \bigoplus_{i=1}^r \Sym^{d_i}S_j^\ast \to \grass(j,N),
\end{equation}
a vector bundle of rank $\nu_j(\bd)\coloneqq\sum_{i=1}^r\binom{d_i+j}{j}$, and $\mcf(j,X)$ is the vanishing locus of this section.

Our goal in this section is to build a bridge from this contemporary perspective to the classical perspective of polars and polar cones which plays a central role in \cite{Sutherland2021}. Notably, this latter perspective appears to be lacking in treatments of $\mcf(j,X)$ in the last several decades, which, as observed in \cite[p. 2]{CilibertoZaidenberg2020} are largely centered on Schubert calculus, a ``trick due to Debarre--Manivel'' \cite{DebarreManivel1998}, and the Bott residue formula and localization in equivariant Chow rings \cite{Hiep2016}. By adding polar cones to the mix, we can get finer information on the arithmetic of $\mcf(j,X)$ over non-algebraically closed $k$, see Remark~\ref{remark:compare_to_DM}.

To study the arithmetic of $\mcf(j,X)$, we introduce an incidence variety parameterizing $j$-planes along with ordered spanning sets.

\begin{definition}
    Let $X\subset\Pb^N$ be a subscheme defined by equations of degree $\bd$. Let $j\ge 0$, and let $J$ be a set of cardinality $j+1$. The incidence variety of spanning sets for $j$-planes in $X$, denoted $\mcc^J(X)$, is defined as
    \[
        \mcc^J(X):=\{((x_a)_{a\in J},\Lambda)\in X^{J}\times \mcf(j,X)~|~\Lambda=\Lambda((x_a)_{a\in J})\}.
    \]
\end{definition}

\begin{remark}
The assignment $J\mapsto \mcc^J(X)$ is contravariant for inclusions $\iota: I\into J$ via the forgetful map
\begin{align*}
    \mcc^J(X)&\to \mcc^I(X)\\
    ((x_a)_{a\in J},\Lambda)&\mapsto ((x_{\iota(b)})_{b\in I},\Lambda( (x_{\iota(b)})_{b\in I}) ).
\end{align*}
In particular, $S_J:=\Perm(J)$ acts on $\mcc^J(X)$ and more generally, $\mcc^\bullet(X)$ forms an $FI^{op}$-variety in the sense of \cite{ChurchEllenbergFarb2015}.
\end{remark}

For our purposes, it suffices to restrict to the standard index sets. We write
\[
    \mcc^j(X):=\mcc^{\{0,\ldots,j\}}(X),
\]
so that $\mcc^0(X)=X$. The inclusion $\{0,\ldots,j-1\}\subset \{0,\ldots,j\}$ defines a map
\[
    \mcc^j(X)\to^{p_j} \mcc^{j-1}(X)
\]
which sits in a natural span
\begin{equation}\label{e:keyspan}
\begin{tikzcd}
    & \mcc^j(X) \arrow[dl, "p_j"'] \arrow[dr, "q_j"] & \\
    \mcc^{j-1}(X) & & \mcf(j,X)
\end{tikzcd}
\end{equation}
where the projection $q_j$ is the map given by $(x_0,\ldots,x_j,\Lambda)\mapsto \Lambda$. This span~\eqref{e:keyspan} will be our key tool to study the arithmetic of $\mcc^j(X)$. Further, for $i\le j$, we denote by
\[
    p_{j,i}: \mcc^j(X)\to \mcc^i(X)
\]
the composite $p_{i+1}\circ\cdots\circ p_j$.

We begin with an explicit description of the fibers of both maps.
\begin{prop}\label{prop:span_properties}
In the span~\eqref{e:keyspan} above,
\begin{enumerate}[label=(\alph*)]
    \item\label{i:qbundle} The map $q_j: \mcc^j(X) \to \mcf(j,X)$ is a Zariski locally trivial fiber bundle with fiber
\[
    q_j^{-1}(\Lambda)=\PGL(\Lambda)/T(\Lambda),
\]
where $T(\Lambda)$ is the stabilizer of an ordered spanning tuple $(x_0,\ldots,x_j)\subset \Lambda$, i.e., the image in $\PGL(\Lambda)\cong\PGL_{j+1}$ of the diagonal torus $\Gb_m^{j+1}$. In particular the fiber is a geometrically irreducible rational variety of dimension $j^2+j$, namely the open subset of $(\Pb^j)^{j+1}$ of spanning tuples. 
    \item\label{i:pfiber} The fiber of $p_j$ at a point $(x_0,\ldots,x_{j-1},\Lambda)\in \mcc^{j-1}(X)$ is given by
\[
    p_j^{-1}(x_0,\ldots,x_{j-1},\Lambda)=S_j|_{\mcf_\Lambda(j,X)}\setminus \Lambda,
\]
i.e. by the complement of the trivial sub-bundle $\Lambda$ in the restriction of the tautological bundle $S_j\to \grass(j,N)$ to $\mcf_\Lambda(j,X)$.
\end{enumerate}
\end{prop}

\begin{proof}
    The first claim follows from identifying $q_j$ with the restriction of the associated bundle of the tautological bundle $S_j\to \grass(j,N)$ to $\mcf(j,X)$. Indeed, let $\mct_{S_j}\to \grass(j,N)$ denote the $\PGL_{j+1}$-torsor underlying $S_j\to \grass(j,N)$. By construction, this torsor is Zariski locally trivial.\footnote{Alternatively, it is the projectivization of a $\GL_{j+1}$-torsor, and thus is Zariski locally trivial by Hilbert's Theorem 90.} Using the fact that $\PGL_{j+1}$ acts simply transitively on ordered tuples of $j+2$ points in general position in $\Pb^j$, i.e., via the identification of schemes
    \[
        \PGL_{j+1}\cong (\Pb^j)^{j+2}-\Sigma,
    \]
    we can identify $\mct_{S_j}$ with the incidence variety
    \[
        \mct_{S_j}=\{(x_0,\ldots,x_{j+1},\Lambda)\in ((\Pb^j)^{j+2}-\Sigma)\times \grass(j,N)~|~x_i\in \Lambda\}.
    \]
    Under this identification, writing $T\cong\Gb_m^j$ for the image of the diagonal torus in $\PGL_{j+1}$, the associated bundle $\mct_{S_j}\times_{\PGL_{j+1}}\PGL_{j+1}/T$ becomes
    \[
        \mct_{S_j}\times_{\PGL_{j+1}}\PGL_{j+1}/T=\{(x_0,\ldots,x_j,\Lambda)\in (\Pb^j)^{j+1}\times \grass(j,N)~|~\Lambda=\Lambda( x_0,\ldots,x_j)\}.
    \]
    Restricting this bundle to $\mcf(j,X)\subset \grass(j,N)$ yields $\mcc^j(X)$, establishing the first claim.

    For the second claim, under the above identification,
    \[
        p_j^{-1}(x_0,\ldots,x_{j-1},\Lambda)=\{(x_0,\ldots,x_{j-1},x_j,\Lambda')\in X^{j+1}\times \mcf(j,X)~|~\Lambda=\Lambda( x_0,\ldots,x_{j-1} ),~\Lambda'=\Lambda( x_0,\ldots,x_j )\}.
    \]
    In particular, $x_j\in \Lambda'\setminus \Lambda$ without any further condition. This implies that the right-hand side is identical to
    \[
        S_j|_{\mcf_\Lambda(j,X)}\setminus \Lambda,
    \]
    and the second claim follows.
\end{proof}

The main result of this section is the following.
\begin{theorem}\label{thm:type_of_fano}
Let $\bd\in \Md$ be a tuple of positive integers with underlying type $\bm$. Let $k$ be a field.
Let $X\subset \Pb^N$ be a subscheme defined over $k$ by equations of degree $\bd$, let $L/k$ be a field extension, and let $\Lambda\subset X_L$ be a $(j-1)$-plane defined over $L$. Let $\bm^j$ be as in Definition~\ref{d:type}. Then:
\begin{enumerate}[label=(\alph*)]
    \item\label{i:tf1} $(S_j|_{\mcf_\Lambda(j,X)}\setminus \Lambda)^{\red}$ is isomorphic to a Zariski open of a reduced subvariety of $\Pb^N_L$ underlying a subscheme defined over $L$ by equations of degree $\bd^j$ and type $\bm^j$, and is dense when non-empty. In particular, $p_j$ is surjective when $N\ge |\bm^j|+j$.
    \item\label{i:tf2} $\mcf_\Lambda(j,X)^{\red}$ is isomorphic to the reduced subvariety of $\Pb^{N-j}_L$ underlying a subscheme defined over $L$ by equations of degree $\bd^j$ and type $\bm^j$. In particular, $\mcf_\Lambda(j,X)$ is nonempty when $N\ge|\bm^j|+j$.
\end{enumerate}
\end{theorem}

We defer the proof for the moment, and record the following corollary.
\begin{corollary}\label{c:tf3}
    We continue the notation of the previous theorem. If $N\geq |\bm^j|+j$, then every point of $X$ (over any extension of $k$) lies on a $j$-plane contained in $X$, i.e., $p_{j,0}: \mcc^j(X)\to X$ is surjective.
\end{corollary}
\begin{proof}
We induct on $j$, the case $j=0$ being trivial. For the induction step, by the inductive hypothesis and the fact that surjectivity is preserved under composition, it suffices to show that $p_j: \mcc^j(X)\to \mcc^{j-1}(X)$ is surjective. By Proposition~\ref{prop:span_properties}~\ref{i:pfiber}, the preimage under $p_j$ of any point in $\mcc^{j-1}(X)$ is precisely the complement of a trivial sub-bundle in the restriction of the tautological bundle $S_j\to \grass(j,N)$ to $\mcf_\Lambda(j,X)$. In particular, the fiber is non-empty if and only if $\mcf_\Lambda(j,X)$ is non-empty. But this follows from Theorem~\ref{thm:type_of_fano}~\ref{i:tf2} by our assumption on $N$ and Lemma~\ref{l:comb}~\ref{i:norm}. 
\end{proof}

\begin{remark}\label{remark:compare_to_DM}
    Debarre and Manivel \cite[Th\'eor\`eme 4.3]{DebarreManivel1998} describe the degree of $\mcf(j,X)$ under the Pl\"ucker embedding of $\grass(j,N)$. This embedding both raises the degree and the ambient dimension substantially, e.g. for a general cubic 3-fold $X$, loc. cit. gives $\mcf(1,X)$ as a degree 45 subvariety of $\Pb^9$ \cite[Fig. 1]{DebarreManivel1998}. In contrast, Theorem~\ref{thm:type_of_fano} shows that for a general $x_0\in X$, we can realize $\mcf_{x_0}(1,X)$ as a degree 6 0-cycle in $\Pb^2$ (and, even better, one arising as an intersection of a conic and cubic). This greater control over the degrees is the principal reason for our interest in the present methods.
\end{remark}

Theorem~\ref{thm:type_of_fano} can be viewed as a restatement of the classical theory of polars and polar cones. For a reader looking for a fuller account, we recommend \cite{Sutherland2021}, \cite[Chapter 1]{Dolgachev2012}, and in particular \cite[Ch.~1, \S Historical Notes]{Dolgachev2012}. To set up the proof of Theorem~\ref{thm:type_of_fano}, we first recall basic details from \cite{Sutherland2021} and establish several key lemmas.

\subsection{Polar cones}\label{subsection:polar_cones}

Let $X=\Vb(f)\subset\Pb^N$ be a hypersurface of degree $d$, and let $x_0 \in \Pb^N$. For any point $y \in \Pb^N$ distinct from $x_0$, we can evaluate $f$ along the line spanned by $x_0$ and $y$ to obtain a homogeneous polynomial in two variables $\lambda$ and $\mu$:
\begin{equation}\label{e:taylor}
    f(\lambda x_0 + \mu y) = \sum_{i=0}^{d} \polar{f}{i}(x_0; y) \lambda^i \mu^{d-i}.
\end{equation}
By construction, each coefficient $\polar{f}{i}(x_0; y)$ is bihomogeneous, with degree $i$ in $x_0$ and degree $d-i$ in $y$. 

\begin{definition}[Polars and Polar Cones]\label{def:polars}\mbox{}
\begin{enumerate}[label=(\alph*)]
    \item Let $f$ be a homogeneous polynomial in $N+1$ variables of degree $d$. For a fixed $x_0\in \Pb^N$, we define the \emph{$i$-th polar of $f$ at $x_0$} to be the homogeneous polynomial $y \mapsto \polar{f}{i}(x_0; y)$ of degree $d-i$, as in~\eqref{e:taylor}.
    \item Let $X=\Vb(f)\subset\Pb^N$ be a hypersurface of degree $\bd = (d)$. For $x_0 \in X$, the \emph{polar cone of $X$ at $x_0$}, denoted $C^1(X;x_0)$, is the subscheme cut out by the vanishing of all polars of positive degree at $x_0$:
    \[
        C^1(X;x_0) := \Vb\left(\polar{f}{0}(x_0; y), \polar{f}{1}(x_0; y), \dots, \polar{f}{d-1}(x_0; y)\right) \subset \Pb^N.
    \]
    The degree of this system of equations is $\bd^1 = (1,2,\dots,d)$ with corresponding type
    \[
        (\underbrace{0,\dots,1}_d)^1 = (\underbrace{1,\ldots,1}_d).
    \]
    \item More generally, if $X \subset \Pb^N$ is an intersection of hypersurfaces defined by a tuple of homogeneous polynomials $\bf = (f_1, \dots, f_r)$, the polar cone is the intersection
    \[ C^1(X;x_0) \coloneqq \bigcap_{i=1}^r C^1(\Vb(f_i);x_0). \]
\end{enumerate}
\end{definition}

\begin{example}
Let $f$ be homogeneous of degree $d$. Then:
\begin{itemize}
\item The $d$-th polar at $x_0$ is the constant $f(x_0)$.
\item The $(d-1)$-st polar at $x_0$ defines the tangent hyperplane at $x_0$.
\item The $0$-th polar at $x_0$ recovers the original polynomial $f(y)$.
\end{itemize}
\end{example}

Because a linear subspace lies on an intersection of varieties exactly when it lies on each component, for a scheme $X\subset \Pb^N$ defined by equations $\bf$ and a linear subspace $\Lambda\subset X$, we have
\begin{align*}
    \mcf(j,X)&=\bigcap_i \mcf(j,\Vb(f_i))\\
    \mcf_{\Lambda}(j,X)&=\bigcap_i \mcf_\Lambda(j,\Vb(f_i)).
\end{align*}

While it appears that this definition relies on the choice of equations defining $X$, we can give a purely intrinsic description. Bertini's Lemma \cite{Bertini1923,Segre1945,Sutherland2021} states that $C^1(X;x_0) \subseteq X$ is a cone with cone point $x_0$, and, hence, for every $x_1 \in C^1(X;x_0) \setminus \{x_0\}$ the line $\Lambda(x_0,x_1)$ lies in $X$; indeed, if all polars $\polar{f}{i}(x_0;x_1)$, $i<d$, vanish and $f(x_0)=0$, then by \eqref{e:taylor} $f$ vanishes identically along $\Lambda(x_0,x_1)$. In particular, for any hyperplane $H\subset \Pb^N\setminus \{x_0\}$, the projection from $x_0$ gives an embedding
\begin{equation}\label{e:planecut}
    C^1(X;x_0)^{\red} \cap H\into \mcf_{x_0}(1,X).
\end{equation}
One of our first results will show that this is in fact an isomorphism of underlying reduced schemes.

\begin{prop}\label{prop:line_iff_point_in_cone}
Let $X \subseteq \Pb^N$ be defined by equations $\bf$ and let $x_0\in X$. Then $C^1(X;x_0)^{\red}$ is the union of all lines in $X$ through $x_0$. In particular,
\begin{enumerate}[label=(\alph*)]
    \item $C^1(X;x_0)^{\red}\subset \Pb^N$ depends only on $X$, not on the choice of defining equations $\bf$.
    \item In the context of the span~\eqref{e:keyspan} for $j=1$, we have:
    \begin{enumerate}[label=(\roman*)]
        \item $p_1^{-1}(x_0)=S_1|_{\mcf_{x_0}(1,X)}\setminus\{x_0\}=C^1(X;x_0)^{\red}\setminus \{x_0\}$,
        \item $q_1(p_1^{-1}(x_0))=\mcf_{x_0}(1,X)$, and
        \item for any hyperplane $H\subset \Pb^N\setminus \{x_0\}$, the embedding~\eqref{e:planecut} is an isomorphism of underlying reduced schemes.
    \end{enumerate}
\end{enumerate}
\end{prop}
\begin{proof}
    The key statement to prove is that $C^1(X;x_0)^{\red}$ is the union of all lines in $X$ through $x_0$; the remaining statements follow directly from definitions. For this, it is enough to consider the case when $X$ is a hypersurface defined by a polynomial $f$, since lines are contained in this intersection if and only if they are contained in each hypersurface individually.

    Let $y\neq x_0$. The line $\Lambda(x_0,y)$ lies in $X$ if and only if $f(\lambda x_0+\mu y)$ vanishes identically, i.e., by \eqref{e:taylor}, if and only if $\polar{f}{i}(x_0;y)=0$ for all $0\leq i\leq d$. Since $\polar{f}{d}(x_0;y)=f(x_0)=0$, this is exactly the condition $y\in C^1(X;x_0)$ of Definition~\ref{def:polars}.
\end{proof}

\begin{prop}\label{prop:type_of_C1}
Let $X\subset \Pb^N$ be a subscheme defined by equations of degree $\bd$.  Let $\bm$ be the type of $\bd$, and let $\bd^1$ and $\bm^1$ be as in Definition~\ref{d:type}.  For any extension $L/k$ and any $x_0 \in X(L)$,
\begin{enumerate}[label=(\alph*)]
    \item $C^1(X;x_0)\subset \Pb^N_L$ is a subscheme defined by equations of degree $\bd^1$ and type $\bm^1$. In particular, $C^1(X;x_0)\setminus \{x_0\}$ is non-empty when $N\ge |\bm^1|+1$.
    \item $\mcf_{x_0}(1,X)^{\red}$ is isomorphic over $L$ to a variety in $\Pb^{N-1}_L$ underlying a closed subscheme defined by equations of type $\bm^1$. In particular, it is non-empty when $N\ge |\bm^1|+1$.
\end{enumerate}
\end{prop}
\begin{proof}
    As observed in Example~\ref{ex:plus}, if $Y$ and $Z$ are subschemes defined by equations of degrees $\bd$ and $\bd'$ and types $\bm$ and $\bm'$, then $Y\cap Z$ is defined by equations of degree $\bd\oplus\bd'$ and type $\bm+\bm'$.  Moreover, by Lemma~\ref{l:comb}, we have
    \[
        (\bm+\bm')^1=\bm^1+\bm'^1.
    \]
    Therefore, just as in the proof of Proposition~\ref{prop:line_iff_point_in_cone}, it is enough to consider the case when $X$ is a hypersurface defined by a degree $d$ polynomial.

    For $x_0\in X(L)$, Definition~\ref{def:polars} shows that $C^1(X;x_0)\subset \Pb^N_L$ is defined by equations of degree $\bd^1$ and type $\bm^1$, as claimed. If $N\ge |\bm^1|+1$, then $\dim C^1(X;x_0)\ge 1$ and thus $C^1(X;x_0)\setminus\{x_0\}\neq\emptyset$.

    The second statement now follows from the isomorphism~\eqref{e:planecut}, as established in Proposition~\ref{prop:line_iff_point_in_cone}.
\end{proof}

\subsubsection*{Iterated polar cones}

Following \cite{Sutherland2021}, we may now iterate this construction. A \emph{$(j-1)$-polar point} of $X$ defined over $L/k$ is a tuple $(x_0,\dots,x_{j-1})\in X^{j}(L)$ of points such that $\Lambda(x_0,\dots,x_{j-1})$ is a $(j-1)$-plane contained in $X$; thus $(j-1)$-polar points are exactly the $L$-points of $\mcc^{j-1}(X)$.\footnote{Note that the difference between our presentation and Sutherland's can essentially be phrased as the difference between a basis and the complete flag associated to it, where we emphasize the former and Sutherland the latter. The emphasis on flags hews closely to the iterative construction of (iterated) polar cones, whereas ours emphasizes the underlying linear algebra. The reader is urged to compare \cite[Definition 2.22]{Sutherland2021} with Proposition \ref{prop:j-polar_points}.}

Rather than iterating \eqref{e:taylor}, it is convenient to expand along a whole linear subspace at once. For a homogeneous $f$ of degree $d$, linearly independent $x_0,\dots,x_{j-1}\in\Pb^N$, and $y\in\Pb^N$, write
\begin{equation}\label{e:taylor_j}
    f(\lambda_0 x_0+\cdots+\lambda_{j-1}x_{j-1}+\mu y) = \sum_{\substack{I\in\Nb^j\\ |I|\leq d}} \polar{f}{I}(x_0,\dots,x_{j-1}; y)\,\lambda^I \mu^{d-|I|},
\end{equation}
where $\lambda^I=\lambda_0^{I_0}\cdots\lambda_{j-1}^{I_{j-1}}$. The coefficient $\polar{f}{I}(x_0,\dots,x_{j-1};y)$, which we call the \emph{$I$-th polar of $f$ at $(x_0,\dots,x_{j-1})$}, is homogeneous of degree $d-|I|$ in $y$ and of multidegree $I$ in $(x_0,\dots,x_{j-1})$. 

\begin{definition}\label{def:iterated_polar_cone}
Let $X\subset\Pb^N$ be defined by $\bf=(f_1,\dots,f_r)$ of degrees $\bd$, and let $(x_0,\dots,x_{j-1})$ be a $(j-1)$-polar point of $X$. The \emph{$j$-th polar cone of $X$ at $(x_0,\dots,x_{j-1})$} is
\[
    C^j(X;x_0,\dots,x_{j-1}) \coloneqq \Vb\Big( \polar{f_i}{I}(x_0,\dots,x_{j-1};\,\cdot\,)  \mid 1\leq i\leq r,\ |I|<d_i \Big)\subset\Pb^N .
\]
\end{definition}

Note that the omitted coefficients, those with $|I|=d_i$, are the coefficients of $f_i|_{\Lambda(x_0,\dots,x_{j-1})}$ and vanish identically precisely because $(x_0,\dots,x_{j-1})$ is a polar point. For $j=1$, the above recovers Definition~\ref{def:polars}. The relation to Sutherland's recursive definition is the following.

\begin{lemma}\label{lemma:iterated_cone_recursion}
For $j\geq 2$ and $(x_0,\dots,x_{j-1})$ a $(j-1)$-polar point of $X$,
\[
    C^j(X;x_0,\dots,x_{j-1}) = C^1\big(C^{j-1}(X;x_0,\dots,x_{j-2});\,x_{j-1}\big)
\]
as subschemes of $\Pb^N$. In particular, iterated polar cones are nested:
\[
    C^{j}(X;x_0,\dots,x_{j-1}) \subseteq C^{j-1}(X;x_0,\dots,x_{j-2}) \subseteq \cdots \subseteq C^1(X;x_0) \subseteq X.
\]
\end{lemma}
\begin{proof}
It suffices to consider a single $f$ of degree $d$. Write $I=(I',a)$ with $I'\in\Nb^{j-1}$, $a\in\Nb$. Substituting $y' = \lambda_{j-1}x_{j-1}+\mu y$ into the expansion \eqref{e:taylor_j} for $(x_0,\dots,x_{j-2})$ and then expanding each $\polar{f}{I'}(x_0,\dots,x_{j-2};\,\cdot\,)$ along the line through $x_{j-1}$ and $y$ via \eqref{e:taylor}, a comparison of coefficients gives
\[
    \polar{f}{(I',a)}(x_0,\dots,x_{j-1};y) = \polar{\big(\polar{f}{I'}(x_0,\dots,x_{j-2};\,\cdot\,)\big)}{\,a}(x_{j-1};y).
\]
As $a$ ranges over $0,\dots,d-|I'|-1$, the right-hand side ranges over all polars of positive degree at $x_{j-1}$ of the form $\polar{f}{I'}(x_0,\dots,x_{j-2};\,\cdot\,)$ of degree $d-|I'|$. Those $\polar{f}{I'}$ with $|I'|<d$ are exactly the equations of $C^{j-1}(X;x_0,\dots,x_{j-2})$, while those with $|I'|=d$ are identically zero. The equality of the two sets of equations follows, and nestedness follows.
\end{proof}

Just as for polar cones, we have the following result for iterated polar cones.

\begin{prop}\label{prop:j-polar_points}
Let $X \subseteq \Pb^N$ be defined by equations $\bf$ and let $(x_0,\dots,x_{j-1})$ be a $(j-1)$-polar point of $X$ defined over $L/k$, spanning $\Lambda=\Lambda(x_0,\dots,x_{j-1})$. Then $C^{j}(X;x_0,\ldots,x_{j-1})^{\red}$ is the union of all $j$-planes in $X$ containing $\Lambda$. In particular,
\begin{enumerate}[label=(\alph*)]
\item The variety $C^{j}(X;x_0,\ldots,x_{j-1})^{\red}\subset \Pb^N$ depends only on $X$ and $\Lambda$, not on the defining forms $\bf$ nor on the chosen spanning points.
\item \label{i:jpolarspan} In the context of the span~\eqref{e:keyspan}, we have
    \begin{enumerate}[label=(\roman*)]
    \item\label{i:ft1a} $p_j^{-1}(x_0,\ldots,x_{j-1},\Lambda)=S_{j}|_{\mcf_{\Lambda}(j,X)}\setminus\Lambda=C^{j}(X;x_0,\ldots,x_{j-1})^{\red}\setminus \Lambda$,
    \item $q_j(p_j^{-1}(x_0,\ldots,x_{j-1},\Lambda))=\mcf_{\Lambda}(j,X)$, and
    \item\label{i:cutisom} for any $(N-j)$-plane $H\subset \Pb^N$ disjoint from $\Lambda$, the map $y\mapsto \langle\Lambda,y\rangle$ gives an isomorphism
    \[
        (C^{j}(X;x_0,\ldots,x_{j-1})\cap H)^{\red} \to^\cong \mcf_{\Lambda}(j,X)^{\red}.
    \]
    \end{enumerate}
\end{enumerate}
\end{prop}
\begin{proof}
As with Proposition~\ref{prop:line_iff_point_in_cone}, it suffices to prove the first claim for a hypersurface $X=\Vb(f)$. If $y\in\Lambda$, then $f(\sum_l\lambda_lx_l+\mu y)$ vanishes identically, so all coefficients in \eqref{e:taylor_j} vanish and $\Lambda\subseteq C^j(X;x_0,\dots,x_{j-1})$. If $y\notin\Lambda$, then $\langle\Lambda,y\rangle=\Lambda(x_0,\dots,x_{j-1},y)$ is a $j$-plane, and it lies in $X$ if and only if $f(\sum_l\lambda_lx_l+\mu y)$ vanishes identically, i.e., if and only if $\polar{f}{I}(x_0,\dots,x_{j-1};y)=0$ for all $|I|\leq d$. Since the coefficients with $|I|=d$ vanish automatically, this is exactly the condition $y\in C^j(X;x_0,\dots,x_{j-1})$. Hence $C^j(X;x_0,\dots,x_{j-1})^{\red}$ is the union of $\Lambda$ and of all $j$-planes of $X$ containing $\Lambda$; and $\Lambda$ is contained in each such $j$-plane (or, if there are none, we simply have $C^j(X;x_0,\dots,x_{j-1})^{\red}=\Lambda$). The remaining statements follow by unpacking the definitions and Proposition~\ref{prop:span_properties}\ref{i:pfiber}.
\end{proof}

Using Lemma~\ref{l:comb}, the multi-index expansion \eqref{e:taylor_j} yields the following directly.\footnote{Note that this generalizes \cite[Proposition 2.26]{Sutherland2021}, which concerns an intersection of type $(1,\dots,1)$.}

\begin{prop}\label{prop:type_of_Cj}
Let $X\subset \Pb^N$ be defined by equations of degree $\bd$.  Let $\bm = t(\bd)$ be the type of $\bd$, and let $\bd^j$ and $\bm^j$ be as in Definition~\ref{d:type}.  For any $L/k$ and any $(j-1)$-polar point $(x_0,\ldots,x_{j-1})$ of $X$ over $L$,
\begin{enumerate}[label=(\alph*)]
    \item\label{i:ft1b} $C^j(X;x_0,\ldots,x_{j-1})\subset \Pb^N_L$ is a subscheme defined by equations of degree $\bd^j$ (up to reordering) and type $\bm^j$. In particular, $C^j(X;x_0,\ldots,x_{j-1})\setminus\Lambda(x_0,\ldots,x_{j-1})$ is non-empty when $N\ge |\bm^j|+j$.
    \item\label{i:ft2} $\mcf_{\Lambda(x_0,\ldots,x_{j-1})}(j,X)^{\red}$ is isomorphic over $L$ to a variety in $\Pb^{N-j}_L$ underlying a closed subscheme defined by equations of degree $\bd^j$ and type $\bm^j$. In particular, it is non-empty when $N\ge |\bm^j|+j$.
\end{enumerate}
\end{prop}
\begin{proof}
    By Lemma~\ref{l:comb}\ref{i:monoid_endo} and Example~\ref{ex:plus} it suffices to treat a single form $f$ of degree $d$, so $\bm=(0,\dots,0,1)$. By Definition~\ref{def:iterated_polar_cone}, the equations of $C^j(X;x_0,\dots,x_{j-1})$ are the polars $\polar{f}{I}$ with $I\in\Nb^j$, $|I|<d$, and $\polar{f}{I}$ has degree $e=d-|I|\in\{1,\dots,d\}$. For fixed $e$, the number of such $I$ is $\#\{I\in\Nb^j: |I|=d-e\}=\binom{d-e+j-1}{j-1}$, which by Lemma~\ref{l:comb}\ref{i:entries} is the $e$-th entry of $\bm^j$. Hence the type is $\bm^j$, and by Lemma~\ref{l:comb}\ref{i:type_commutes} the degrees are those of $\bd^j$ up to order. The dimension statement is Krull's theorem: 
    \[ 
        \dim C^j(X;x_0,\dots,x_{j-1})\geq N-|\bm^j|\geq j 
    \]
    forces the cone to be strictly larger than the $(j-1)$-plane $\Lambda(x_0,\dots,x_{j-1})$. The second part now follows from Proposition~\ref{prop:j-polar_points}\ref{i:jpolarspan}\ref{i:cutisom}, since $H\cong\Pb^{N-j}_L$ and $C^j\cap H$ is cut out by the restrictions of the same equations. 
\end{proof}

\begin{proof}[Proof of Theorem~\ref{thm:type_of_fano}]
    After extending scalars to $L$ and choosing a spanning tuple $(x_0,\dots,x_{j-1})$ of $\Lambda$ over $L$, the first statement follows by combining statement~\ref{i:jpolarspan}\ref{i:ft1a} of Proposition~\ref{prop:j-polar_points} with statement~\ref{i:ft1b} of Proposition~\ref{prop:type_of_Cj}. The second statement of the theorem is statement~\ref{i:ft2} of Proposition~\ref{prop:type_of_Cj}.
\end{proof}

\section{Parameter spaces and general intersections}\label{s:param}

The results of Section~\ref{s:geom} concern a fixed intersection $X$. The inductive arguments of Section~\ref{s:solvable} will, however, only produce dense rational points on \emph{general} intersections of a given type, and the induction must therefore be carried out at the level of the parameter space of all intersections of that type. In this section we set up that parameter space, the groups acting on it, and the notion of ``general'' we use. We also prove the two results permitting induction over parameter spaces: the identification of the polar model of a general intersection containing a fixed plane with a general intersection in a smaller projective space (Theorem~\ref{thm:linear_model}), and a bound on the dimension of the fibers of a family over a general point (Lemma~\ref{lemma:general-dim}). 

\subsection{Universal families}\label{subsection:universal}

For $N \geq 0$ and $\bd = (d_1, \dots, d_r) \in \Md$ with each $d_i \geq 1$, we write
\[
  \param{N}{\bd} \coloneqq \prod_{i=1}^r H^0\big(\Pb^N_k, \mco(d_i)\big)
\]
for the affine space over $k$ of tuples $\bf = (f_1, \dots, f_r)$ of forms of degrees $\bd$ in $z_0,\dots,z_N$. By convention, we take $\param{N}{\mathbf{0}}$ to be a point. For a field $L/k$ and $\bf \in \param{N}{\bd}(L)$, we write $X_{\bf} \coloneqq \Vb(f_1,\dots,f_r) \subseteq \Pb^N_L$ and denote the type by $\bm=t(\bd)$. Note that $\GL_{N+1}$ acts on $\param{N}{\bd}$ by substitution:
\[
  (g \cdot \bf)(z) \coloneqq \bf(g^{-1} z), \qquad\text{so that}\qquad X_{g \cdot \bf} = g (X_\bf).
\]

The \emph{universal intersection of degree $\bd$} is the closed subscheme
\[
  \mcx = \mcx_N(\bd) \coloneqq \{(\bf,x)\in\param{N}{\bd}\times\Pb^N \mid x\in X_\bf\},
\]
whose fiber over $\bf$ is $X_\bf$. More generally, for $j\geq0$ we have the \emph{relative scheme of $j$-planes}
\[
  \mcf(j,\mcx) \coloneqq \{(\bf,\Lambda)\in\param{N}{\bd}\times\grass(j,N) \mid \Lambda\subseteq X_\bf\},
\]
so that $\mcf(0,\mcx)=\mcx$, the fiber of $\mcf(j,\mcx)\to\param{N}{\bd}$ over $\bf$ is $\mcf(j,X_\bf)$, and $\mcf(j,\mcx)$ is the incidence variety denoted $I_j$ in \cite[\S2]{DebarreManivel1998}. We have the incidence correspondence
\begin{equation}\label{eqn:universal_incidence}
\begin{tikzcd}
& \mcf(j,\mcx) \arrow[swap]{dl} \arrow{dr} \\
\param{N}{\bd} && \grass(j,N).
\end{tikzcd} 
\end{equation}
Containing a given $j$-plane is a linear condition on the coefficients of $\bf$. In particular, we will use that a degree $d$ form $f$ vanishes on $\Lambda_0 \coloneqq \Vb(z_{j+1},\dots,z_N) \in \grass(j,N)$ if and only if the $\binom{d+j}{j}$ coefficients of degree $d$ monomials purely in $z_0,\dots,z_j$ vanish. Hence the fiber over $\Lambda_0$,
\begin{equation}
    \param{N, \Lambda_0}{\bd} \coloneqq \{\bf \in \param{N}{\bd} \mid \Lambda_0 \subseteq X_\bf\},
\end{equation}
is the subspace of $\param{N}{\bd}$ cut out by the vanishing of these coefficients for each $f_i$, of codimension
\[
  \nu_j(\bd)\coloneqq\sum_{i=1}^r\binom{d_i+j}{j}.
\]

\begin{lemma}\label{lemma:universal_smooth}
The projection $\mcf(j,\mcx)\to\grass(j,N)$ is Zariski locally trivial. Thus $\mcf(j,\mcx)$ is smooth and irreducible of dimension $\dim\param{N}{\bd} + \delta_j$, where
\begin{equation}\label{eqn:expected_dimension}
    \delta_j = \delta(N,\bd,j) \coloneqq (j+1)(N-j) - \nu_j(\bd).
\end{equation}
Moreover, for every $\bf \in \param{N}{\bd}$, every irreducible component of $\mcf(j,X_\bf)$ has dimension at least $\delta_j$; if $\delta_j < 0$ then a general $X_\bf$ contains no $j$-planes. 
\end{lemma}
\begin{proof}
The projection $\mcf(j,\mcx)\to\grass(j,N)$ is $\GL_{N+1}$-equivariant. Moreover, $\GL_{N+1}$ acts transitively on $\grass(j,N)$ and the orbit map $g \mapsto g \cdot \Lambda_0$ has a section $\Lambda \mapsto g_{\Lambda}$ over each standard affine chart $V \subseteq \grass(j,N)$. Since $\Lambda \subseteq X_\bf$ if and only if $\Lambda_0 \subseteq X_{g_{\Lambda}^{-1} \cdot \bf}$, the automorphism of $\param{N}{\bd}\times V$ given by
\begin{equation}\label{eqn:trivialization}
  (\bf, \Lambda) \mapsto (g_{\Lambda}^{-1} \cdot \bf, \Lambda)
\end{equation}
carries $\mcf(j,\mcx)|_V$ onto $\param{N,\Lambda_0}{\bd} \times V$. This gives local triviality, and that $\mcf(j,\mcx)$ is smooth of pure dimension 
\[
  \dim \param{N}{\bd} - \nu_j(\bd) + (j+1)(N-j) = \dim\param{N}{\bd}+\delta_j.
\] 
Lastly, since $\mcf(j,X_\bf)$ is the fiber of $\mcf(j,\mcx) \to \param{N}{\bd}$ over $\bf$, it is cut out of $\mcf(j,\mcx)$ by the $\dim\param{N}{\bd}$ linear equations defining $\bf$. By Krull's theorem, each of its irreducible components has dimension at least
\[ 
    \dim\mcf(j,\mcx)-\dim\param{N}{\bd} = \delta_j.
\]
If $\delta_j < 0$, then $\dim\mcf(j,\mcx)<\dim\param{N}{\bd}$, so the projection is not dominant and $\mcf(j,X_\bf)=\emptyset$ for $\bf$ outside the closure of its image.
\end{proof}

\subsection{Groups and general intersections}\label{subsection:general}
Fix $N \geq 0$ and $\bd = (d_1, \dots, d_r) \in \Md$ with each $d_i \geq 1$. By reordering the polynomials if necessary, we assume throughout that $d_1 \leq \cdots \leq d_r$.

In addition to the general linear group $\GL_{N+1}$ acting by substitution of coordinates, we will make use of the algebraic group $\mcr = \mcr_N(\bd)$ of \emph{changes of generators}. $\mcr$ consists of the invertible $r \times r$ matrices $\theta = (\theta_{ii'})$ with each $\theta_{ii'} \in k[z_0, \dots, z_N]$ homogeneous of degree $d_i - d_{i'}$ (so $\theta_{ii'}=0$ if $d_i<d_{i'}$). $\mcr$ acts on $\param{N}{\bd}$ by
\[
  (\theta\cdot\bf)_i \coloneqq \sum_{i'=1}^r \theta_{ii'}\,f_{i'} .
\]
Since $d_1\leq\cdots\leq d_r$, such a matrix is block lower-triangular with constant diagonal blocks $(\theta_{ii'})_{d_i=d_{i'}=d}$, and it is invertible if and only if these blocks are. Thus, if $\bm=(m_1,\dots,m_n)=t(\bd)$, we obtain a Levi-type decomposition
\[
  \mcr_N(\bd) \cong \Big( \prod_{i=1}^{n} \GL_{m_i} \Big) \ltimes \mcu ,
\]
where the unipotent radical $\mcu$ consists of matrices with identity blocks on the diagonal and homogeneous polynomials of strictly positive degree below the diagonal. Note that $\theta\cdot\bf$ is a new system of generators of the ideal $(f_1,\dots,f_r)$, so that $X_{\theta\cdot\bf}=X_\bf$ as subschemes.
Substitution of coordinates preserves degrees, so $\GL_{N+1}$ normalizes $\mcr$, i.e., $g\theta g^{-1}$ has entries given by $\theta_{ii'}(g^{-1}z)$.

\begin{definition}\label{definition:graded_equivalence}
For $N \geq 0$ and $\bd$ as above, we set
\[
  \mcg = \mcg_N(\bd) \coloneqq \GL_{N+1} \ltimes \mcr_N(\bd),
\]
acting on $\param{N}{\bd}$ by $(g,\theta)\cdot\bf = g\cdot(\theta\cdot\bf)$. 
\end{definition}

Morally, one considers $\mcg$ because one is interested in properties of tuples $\bf$ that depend only on the projective equivalence class of the subscheme $X_{\bf}\subseteq\Pb^N$; any such property is $\mcg$-invariant. Practically, the group $\mcg$ acts on the universal families of Section~\ref{subsection:universal}: on $\mcx$ and on $\mcf(j,\mcx)$, $\GL_{N+1}$ acts diagonally and $\mcr$ acts on the first factor only, and the projections to $\param{N}{\bd}$ are $\mcg$-equivariant. 

\begin{definition}\label{definition:general_property}
Let $K/k$ be an infinite field extension, let $\bd\in\Md$ with $d_i\geq1$, and write $\bm = t(\bd)$. Let $\mcp$ be a property of tuples $\bf \in \param{N}{\bd}(K)$. We say that $\mcp$ \emph{holds for general $\bf$ of type $\bm$}, or that \emph{a general intersection $X\subseteq\Pb^N_K$ of type $\bm$ satisfies $\mcp$}, if there is a dense Zariski open subscheme $U \subseteq \param{N}{\bd}$ (over $k$) such that $\mcp(\bf)$ holds for every $\bf \in U(K)$.
\end{definition}

\begin{remark}\label{remark:general}\mbox{}
\begin{enumerate}
\item\label{i:general_dense} Since $\param{N}{\bd}$ is irreducible, any such $U$ is dense, and finitely many properties that each hold for general $\bf$ hold simultaneously for general $\bf$. 
\item\label{i:general_k} Since $U(k)\subseteq U(K)$, if $\mcp$ holds for general $\bf$ of type $\bm$ then in particular it holds for all $k$-rational tuples $\bf$ outside a proper closed subset of $\param{N}{\bd}$. Note that there are such $\bf$ whenever $k$ is infinite. Requiring the property over $K$ rather than over $k$ is what permits the inductions of Section~\ref{s:solvable}, as the inductive hypotheses will be applied to intersections defined over $K$.
\item\label{i:general_smooth} If $N\geq r$, the locus of $\bf$ for which $X_\bf$ is a smooth complete intersection of codimension $r$ is a dense open subset of $\param{N}{\bd}$; hence in Definition~\ref{definition:general_property} one may always shrink $U$ to lie inside the locus of smooth complete intersections.
\end{enumerate}
\end{remark}

\subsection{The linear model of \texorpdfstring{$\mcf_{\Lambda_0}(j,X)$}{F\_Lambda0(j,X)}}\label{subsection:linear_model} 
For a $(j-1)$-plane $\Lambda\subseteq X$, Theorem \ref{thm:type_of_fano} identifies $\mcf_\Lambda(j,X)^{\red}$ with a subvariety of $\Pb^{N-j}$ cut out by equations of degree $\bd^j$, which we call the \emph{polar model} of $X$ at $\Lambda$. In Section~\ref{s:solvable}, we will apply an inductive hypothesis about general intersections of degree $\bd^j$ in $\Pb^{N-j}$ to polar models, and for this we must know that if $X$ is general among intersections containing $\Lambda$, then its polar model at $\Lambda$ is general too. Here we show that the equations of the polar model are obtained from those of $X$ by a linear isomorphism of parameter spaces, so dense opens correspond to dense opens, and that this isomorphism is compatible with the actions above.

Fix $j \ge 1$ and let $\Lambda_0=\Vb(z_j,\dots,z_N)$ be the standard $(j-1)$-plane, spanned by the coordinate points $e_0,\dots,e_{j-1}$, so that $\param{N,\Lambda_0}{\bd}$ is the fiber of $\mcf(j-1,\mcx)$ over $\Lambda_0$. We write $H_0=\Vb(z_0,\dots,z_{j-1})\cong\Pb^{N-j}$ for the complementary coordinate plane, with coordinates $w = (w_0,\dots,w_{N-j})$; by \eqref{e:taylor_j}, every form $f$ of degree $d$ on $\Pb^N$ is written uniquely as
\begin{equation}\label{eq:expansion}
  f = \sum_{\substack{I \in \Nb^j \\ |I| \le d}} z^I\, \polar{f}{I}(w),
\end{equation}
where $z^I = z_0^{I_0} \cdots z_{j-1}^{I_{j-1}}$ and $\polar{f}{I}(w)$ is the $I$-th polar of $f$ at $(e_0,\dots,e_{j-1})$ restricted to $H_0$.
The constants $\polar{f}{I}$ with $|I|=d$ are the coefficients of the monomials purely in $z_0,\dots,z_{j-1}$, so $\param{N,\Lambda_0}{\bd}$ is the linear subspace cut out by $\polar{f_i}{I} = 0$ for all $|I| = d_i$ as in Section~\ref{subsection:universal}. We write
\[
\begin{aligned}
  \Pi_0: \param{N,\Lambda_0}{\bd} & \to \param{N-j}{\bd^j},\\
  \bf & \mapsto \big(\polar{f_i}{I}\big)_{1 \le i \le r,\ |I| < d_i},
\end{aligned}
\]
where the forms $\polar{f_i}{I}$ are regarded as forms on $H_0 = \Pb^{N-j}$ and listed in a fixed order of weakly increasing degree; their degrees are exactly those of $\bd^j$ by Proposition~\ref{prop:type_of_Cj}. Lastly, we write $\mcg_{\Lambda_0}$ for the stabilizer of $\Lambda_0$, consisting of the pairs $(g,\theta)$ with $g\cdot\Lambda_0=\Lambda_0$.

\begin{theorem}\label{thm:linear_model}\mbox{}
\begin{enumerate}[label=(\alph*)]
\item\label{i:lm_linear} $\Pi_0$ is a $k$-linear isomorphism.
\item\label{i:lm_iso} For every field $L/k$ and $\bf\in\param{N,\Lambda_0}{\bd}(L)$, the map $y\mapsto\langle\Lambda_0,y\rangle$ is an isomorphism of $L$-schemes
\[
  X_{\Pi_0(\bf)}^{\red}\xrightarrow{\ \sim\ }\mcf_{\Lambda_0}(j,X_\bf)^{\red}.
\]
\item\label{i:lm_equivariant} For every field $L/k$ and $h\in\mcg_{\Lambda_0}(L)$, there is $h'\in\mcg_{N-j}(\bd^j)(L)$ with $\Pi_0(h\cdot\bf)=h'\cdot\Pi_0(\bf)$ for all $\bf\in\param{N,\Lambda_0}{\bd}$. Hence, if $U'\subseteq\param{N-j}{\bd^j}$ is a $\mcg_{N-j}(\bd^j)$-stable open, then $\Pi_0^{-1}(U')\subseteq\param{N,\Lambda_0}{\bd}$ is a $\mcg_{\Lambda_0}$-stable open.
\end{enumerate}
\end{theorem}
\begin{proof}
For \ref{i:lm_linear}, the coordinates of $\param{N}{\bd}$ are partitioned into the constants $\polar{f_i}{I}$, $|I|=d_i$, which vanish on $\param{N,\Lambda_0}{\bd}$, and the coefficients of the forms $\polar{f_i}{I}$, $|I|<d_i$, which are the coordinates of $\param{N-j}{\bd^j}$.
For \ref{i:lm_iso}, 
\[ 
  X_{\Pi_0(\bf)}=C^j(X_\bf;e_0,\dots,e_{j-1})\cap H_0
\] 
by the description of $\Pi_0(\bf)$ above: the equations $\Pi_0(\bf)$ are precisely the restrictions to $H_0$ of the polars of $\bf$ whose vanishing defines $C^j(X_{\bf};e_0,\dots,e_{j-1})$, i.e. this is Proposition~\ref{prop:j-polar_points}\ref{i:cutisom} with $H=H_0$.

For \ref{i:lm_equivariant}, recall that $\bf\in\param{N}{\bd}$ is determined by its coefficients $\polar{f_i}{I}$ for $1\leq i\leq r$ and $|I|\leq d_i$, which are forms on $H_0$ of degrees $d_i-|I|$: those with $|I|<d_i$ have their degrees described by $\bd^j$, and those with $|I|=d_i$ have degree $0$. The subspace $\param{N,\Lambda_0}{\bd}$ is where the entries of degree $0$ vanish and $\Pi_0$ forgets them. By~\ref{i:lm_linear}, we see that the entries of $\Pi_0(h\cdot\bf)$ only depend on the entries of $\Pi_0(\bf)$. Since $h\in\mcg_{\Lambda_0}(L)$ fixes $\Lambda_0$, it preserves $\param{N,\Lambda_0}{\bd}$. Therefore, as a matrix in the coordinates $(f_i^{[I]})$ of $\param{N}{\bd}$, the map induced on $\param{N-j}{\bd^j}$ via~\ref{i:lm_linear} is the matrix $h'$ obtained by omitting the rows and columns of the degree-$0$ entries. An elementary check verifies that $h'\in 
\mcg_{N-j}(\bd^j)$ as claimed, and the ``hence'' now follows from the equivariance of $\Pi_0$.
\end{proof}

\subsection{Dimension counts}\label{subsection:genericity}
We begin with the schemes $\mcf(j,X_\bf)$ themselves. Recall from \cite[\S1]{DebarreManivel1998} the quantity defined in terms of \eqref{eqn:expected_dimension}:
\[
  \delta^-(N,\bd,j) \coloneqq \min\{\delta(N,\bd,j),\ N-2j-r\}.
\]

\begin{theorem}[{\cite[Th\'eor\`eme~2.1]{DebarreManivel1998}}]\label{thm:DM}
Let $X\subseteq\Pb^N_{\bar k}$ be a subscheme defined by equations of degree $\bd$. 
\begin{enumerate}[label=(\alph*)]
\item\label{i:DM_empty} If $\delta^-(N,\bd,j) < 0$, then $\mcf(j,X)=\emptyset$ for general $X$.
\item\label{i:DM_smooth} If $\delta^-(N,\bd,j)\geq0$, then $\mcf(j,X)\neq\emptyset$, and $\mcf(j,X)$ is smooth of dimension $\delta(N,\bd,j)$ for general $X$.
\item\label{i:DM_connected} If $\delta^-(N,\bd,j)>0$, then $\mcf(j,X)$ is connected.
\end{enumerate}
\end{theorem}

``General'' refers to a dense open subset of $\param{N}{\bd}_{\bar k}$. We record what we need in the sense of Definition~\ref{definition:general_property}:

\begin{corollary}\label{cor:DM_locus}
Let $\Delta_j\subseteq\param{N}{\bd}$ be the image of the closed subset on which $\pi:\mcf(j,\mcx)\to\param{N}{\bd}$ is not smooth. Then $\Delta_j$ is closed, defined over the prime field, and $\mcg_N(\bd)$-stable. Moreover, for every field $L/k$ and $\bf\in(\param{N}{\bd}\setminus\Delta_j)(L)$, the scheme $\mcf(j,X_\bf)$ is smooth of pure dimension $\delta_j$. If $\delta^-(N,\bd,j)\geq0$ then $\Delta_j$ is proper, so this holds for general $\bf$ of type $\bm$; and if $\delta^-(N,\bd,j)>0$, then $\mcf(j,X_\bf)$ is non-empty and geometrically irreducible for every such $\bf$.
\end{corollary}
\begin{proof}
Let $S\subseteq\mcf(j,\mcx)$ be the closed subset on which $\pi$ is not smooth. By the properness of $\pi$, we know $\Delta_j = \pi(S)$ is closed. Since $\pi$ is defined over the prime field, so are $S$ and $\Delta_j$. Similarly, by equivariance of $\pi$ and $\mcg$-stability of $S$, $\Delta_j$ is $\mcg$-stable.

By the smoothness of $\mcf(j,\mcx)$ and $\param{N}{\bd}$, $\pi$ is smooth at $(\bf,\Lambda)$ if and only if the tangent space to the fiber $\mcf(j,X_\bf)$ at $\Lambda$ has dimension
\[ 
    \dim\mcf(j,\mcx)-\dim\param{N}{\bd} = \delta_j.
\]
Every component of $\mcf(j,X_\bf)$ has dimension at least $\delta_j$ by Lemma~\ref{lemma:universal_smooth}, so this happens if and only if $\mcf(j,X_\bf)$ is smooth of dimension $\delta_j$ at $\Lambda$. Hence $\bf\not\in\Delta_j$ if and only if $\mcf(j,X_\bf)$ is smooth of pure dimension $\delta_j$. 

If $\delta^-(N,\bd,j) \geq 0$, Theorem~\ref{thm:DM}\ref{i:DM_smooth} gives $\bf\in\param{N}{\bd}(\bar k)$ with $\mcf(j,X_\bf)$ smooth of dimension $\delta_j$, and so $\Delta_j$ is proper. Hence $\param{N}{\bd}\setminus\Delta_j$ is a non-empty $\mcg$-stable open as in Definition~\ref{definition:general_property}. Lastly, if $\delta^-(N,\bd,j)>0$ and $\bf\notin\Delta_j$, then $\mcf(j,X_\bf)_{\bar L}$ is non-empty and connected by Theorem~\ref{thm:DM}, and smooth, hence irreducible.
\end{proof}

We next state a lemma allowing us to control the geometry of a general member of a family. 

\begin{lemma}\label{lemma:general-dim}
Let $B$ be an irreducible $k$-variety, $\delta\in\Zb$, and $\pi: E \to B$ a morphism of $k$-varieties with $\dim E \le \dim B + \delta$. Then there is a dense open subset $U \subseteq B$ such that $\dim \pi^{-1}(b) \le \delta$ for every point $b \in U$, with the convention $\dim\emptyset=-\infty$; in particular, if $\delta<0$ then $\pi^{-1}(U)=\emptyset$. Moreover, if a group $\mathbf{G}$ acts on $E$ and $B$ and $\pi$ is $\mathbf{G}$-equivariant, then $U$ may be taken to be $\mathbf{G}$-stable.
\end{lemma}
\begin{proof}
The set $F\coloneqq\{b\in B\mid\dim\pi^{-1}(b)\geq\delta+1\}$ is constructible \cite[\href{https://stacks.math.columbia.edu/tag/05F7}{Lemma 05F7}]{stacks-project}, and it is not dense, since otherwise it would contain a dense open $V$ with $\dim\pi^{-1}(V)\geq\dim B+\delta+1>\dim E$. Take $U\coloneqq B\setminus\overline F$, which is $\mathbf{G}$-stable when $\pi$ is equivariant.
\end{proof}

Special cases of the following appear as \cite[Theorem 2.1]{DebarreManivel1998} and \cite[Theorem 3.2]{cheng2026}. We give the proof for completeness.
\begin{prop}\label{prop:equidim}
Let $j \ge 0$ and $N \ge 0$. For general $\bf \in \param{N}{\bd}$ (over any infinite $K$), every irreducible component of $\mcf(j, X_\bf)$ has dimension exactly $\delta_j=\delta(N,\bd,j)$, with the convention that $\mcf(j, X_\bf) = \emptyset$ if $\delta_j < 0$; and every irreducible component of $\mcc^j(X_\bf)$ has dimension exactly $\delta_j + j^2 + j$.
\end{prop}
\begin{proof}
As noted in Lemma~\ref{lemma:universal_smooth}, every irreducible component of $\mcf(j, X_\bf)$ has dimension at least $\delta_j$ by Krull's theorem. For the upper bound, apply Lemma \ref{lemma:general-dim} to the $\mcg$-equivariant projection $\mcf(j,\mcx) \to \param{N}{\bd}$, whose fiber over $\bf$ is $\mcf(j, X_\bf)$ and whose total space has dimension $\dim\param{N}{\bd}+\delta_j$ by Lemma~\ref{lemma:universal_smooth}. The statement for $\mcc^j(X_\bf)$ follows from Proposition~\ref{prop:span_properties}\ref{i:qbundle}: $q_j:\mcc^j(X_\bf) \to \mcf(j, X_\bf)$ is a Zariski locally trivial bundle with geometrically irreducible fibers of dimension $j^2 + j$, so its irreducible components are the preimages of those of $\mcf(j, X_\bf)$.
\end{proof}

\section{Solvable \texorpdfstring{$j$}{j}-planes via total \texorpdfstring{$j$}{j}-polar spaces}\label{s:solvable}
 
For the remainder of the paper, $k$ is a field with $\characteristic(k)\neq 2,3$ and $K\coloneqq k^{\Sol}$ is a solvable closure. In particular, $K$ is infinite, so all of Section~\ref{s:param} applies, and the results of Section~\ref{s:geom} apply to all intersections of quadrics, cubics and quartics.\footnote{The only other property of $K$ we use is that every one-variable polynomial of degree at most $4$ over $K$ splits into linear factors. In particular, one could weaken ``solvably closed'' to ``closed under quadratic and cubic extensions.''} For $(m_2,m_3,m_4)\in\Nb^3$, we write $\bm=(0,m_2,m_3,m_4)$ for the corresponding type and $\bd$ for its degree tuple, and we use the formulas for $\bm^j$ and $|\bm^j|$ of Example~\ref{example:type_shift}; in particular
\begin{equation}\label{eqn:floor}
    |\bm^j|+j = j+(j+1)m_2+\tbinom{j+2}{2}m_3+\tbinom{j+3}{3}m_4 .
\end{equation}
 
\subsection{The functions \texorpdfstring{$\density{j}{\cdot}$}{f\_j} and \texorpdfstring{$\existence{j}{\cdot}$}{g\_j}}\label{subsection:fj}
 
\begin{definition}\label{definition:fj}
We define $f: \Nb \times \Nb^3 \to \Nb \cup \{\infty\}$ by
\begin{equation}\label{eqn:def_of_fj}
\begin{aligned}
    \density{j}{m_2,m_3,m_4}\coloneqq\min\{ N_0 \mid\ &\text{for every } N \geq N_0, \text{ for general } X \subseteq \Pb^N_{K} \text{ of type } (0,m_2,m_3,m_4), \\
    &\quad \mcf(j,X)(K) \subseteq \mcf(j,X) \text{ is nonempty and dense, and }\\
    &\quad p_{j,0}: \mcc^j(X)\to X \text{ is surjective}\},
\end{aligned}
\end{equation}
with the convention that $\min\emptyset=\infty$.
\end{definition}

Here ``general'' is in the sense of Definition~\ref{definition:general_property}: the two properties in \eqref{eqn:def_of_fj} depend only on the projective equivalence class of $X_\bf\subseteq\Pb^N_K$, hence are $\mcg(K)$-invariant. Since $\mcg(K)\subset\mcg$ is dense, we could, without loss of generality, restate the requirement to hold for all $\bf$ in some non-empty $\mcg$-stable open subset of $\param{N}{\bd}_k$ (since any dense open in which the requirement holds can be replaced with its $\mcg$-orbit). We omit doing so as to avoid cluttering proofs with an extra, unnecessary check.

\begin{remark}\label{remark:other_functions}\mbox{}
\begin{enumerate}
    \item Variations of the function $\density{0}{m_2,m_3,m_4}$, particularly with regards to the \emph{existence} of points rather than their density, have been considered in the literature (e.g., \cite{Sylvester1887,Heberle2021}) for other classes of points; see also Section \ref{subsection:existence}. 
    \item We omit linear forms from the definition of $\density{j}{m_2,m_3,m_4}$ because their inclusion reduces to the case of a smaller ambient space. We remark on this here because the notion of a general intersection involves the whole parameter space, including the linear forms.
    \item Every object of Sections~\ref{s:geom} and \ref{s:param}---the parameter spaces $\param{N}{\bd}$, the universal family and its relative Fano schemes, the groups $\GL_{N+1}$ and $\mcr_N(\bd)$, the polar model, and the linear model $\Pi_0$---is defined over $\Zb$, and the loci $\Delta_j$ of Corollary~\ref{cor:DM_locus} are defined over the prime field. Each open $U\subseteq\param{N}{\bd}$ constructed in this section is cut out from these by conditions that are likewise defined over the prime field, and so the open set in \eqref{eqn:def_of_fj} may be taken to depend on $k$ only through its characteristic. While we have not done so for expositional reasons, one could restate the results of this section over $\spec\Zb[1/6]$, obtaining a $\mcg_N(\bd)$-stable open $U\subseteq\param{N}{\bd}_{\Zb[1/6]}$ such that $X_u$ has a dense set of $\kappa(u)^{\Sol}$-points for every point $u\in U$, of whatever residue field $\kappa(u)$. 
\end{enumerate}
\end{remark}

We record the following lemma for later use. Note that this follows immediately from observing that the non-empty dense open $U\subset \Ab_N(\bd)$ required for Definition~\ref{definition:fj} for a $\bd$ of type $(0,m_2,m_3,m_4)$, pulls back along the ``intersection'' (rational) map 
\[
    \cap\colon\Ab_{N+m_1}((1)^{\oplus m_1})\times\Ab_{N+m_1}(\bd)\dashrightarrow \Ab_{N}(\bd)
\]
to a dense open with the same properties.

\begin{lemma}\label{lemma:omit_linear}
Let $\bm = (m_1, m_2, m_3, m_4) \in \Mt$ with $m_1 \ge 1$, and assume $\density{j}{m_2, m_3, m_4} < \infty$. For every $N \geq \density{j}{m_2, m_3, m_4} + m_1$, a general intersection $X \subseteq \Pb^N_K$ of type $\bm$ satisfies the two properties of \eqref{eqn:def_of_fj}: $\mcf(j,X)(K)\subseteq\mcf(j,X)$ is non-empty and dense, and $p_{j,0}:\mcc^j(X)\to X$ is surjective.
\end{lemma}

We next record some elementary properties of $\density{j}{m_2,m_3,m_4}$.
\begin{lemma}\label{lemma:base_case}
For every $j\in\Nb$, $\density{j}{0,0,0}=j$.
\end{lemma}
\begin{proof}
Here $\bm=\mathbf{0}$ and $\param{N}{\mathbf{0}}$ is a point, so ``general'' is vacuous and we have $X=\Pb^N_K$. For $N<j$, we have $\mcf(j,\Pb^N_K) = \emptyset$, so the first property of \eqref{eqn:def_of_fj} fails and no $N_0<j$ can occur in \eqref{eqn:def_of_fj}. For $N\geq j$, $\mcf(j,\Pb^N_K)=\grass(j,N)_K$ is a non-empty rational variety, whose $K$-points are dense since $K$ is infinite, and $p_{j,0}$ is surjective since every point of $\Pb^N$ lies on a $j$-plane.
\end{proof}
 
\begin{lemma}\label{lemma:nonempty}
For any $m_2,m_3,m_4\in \Nb$, we have $m_2+m_3+m_4 \leq \density{0}{m_2,m_3,m_4}$.
\end{lemma}
\begin{proof}
Set $r \coloneqq m_2+m_3+m_4$ and let $N < r$. The locus $U_\emptyset\subseteq\param{N}{\bd}$ of tuples $\bf$ with $X_\bf=\emptyset$ is open, since its complement is the image of the universal intersection $\mcx$ under the projection $\mcx\to\param{N}{\bd}$. It is also non-empty: for instance, the tuple with $f_i=z_{\min(i-1,N)}^{d_i}$ for $1\leq i\leq r$ lies in $U_\emptyset$, because $r\geq N+1$. If the properties of \eqref{eqn:def_of_fj} held for general $X$ at this $N$, via a non-empty open $U$, then $U\cap U_\emptyset$ would be a non-empty open subset of an affine space over $k$, hence would contain a $K$-point $\bf$ since $K$ is infinite; but $X_\bf(K)=\emptyset$, a contradiction. Hence no $N_0\leq N$ can occur in \eqref{eqn:def_of_fj}; since $N<r$ was arbitrary, $\density{0}{m_2,m_3,m_4}\geq r$.
\end{proof}
 
\begin{lemma}\label{lemma:monotonicity}
Fix $m_2 \leq m'_2$, $m_3 \leq m'_3$, and $m_4 \leq m'_4$. Then
\[
    \density{0}{m_2,m_3,m_4} \leq \density{0}{m'_2,m'_3,m'_4}.
\]
\end{lemma}
\begin{proof}
Let $N\ge \density{0}{m'_2,m'_3,m'_4}$. It suffices to show that $N$ satisfies the conditions for $\density{0}{m_2,m_3,m_4}$ as well. Let $\bd'$ and $\bd$ be the degree tuples corresponding to $\bm'$ and $\bm$. Forgetting the extra polynomials gives a linear projection
\[
    \Ab_N(\bd')\to \Ab_N(\bd)
\]
which fits into a commuting diagram
\[
    \xymatrix{
        \mcx_N(\bd') \ar[r] \ar[d] & \mcx_N(\bd) \ar[d] \ar[r] & \Pb^N \\
        \Ab_N(\bd') \ar[r]^\pi & \Ab_N(\bd) 
    }
\]
From the construction in Section~\ref{subsection:universal}, under the upper horizontal maps to $\Pb^N$, $\mcx_N(\bd')$ and $\mcx_N(\bd)$ are vector sub-bundles of the trivial bundles $\param{N}{\bd'}\times\Pb^N$ and $\param{N}{\bd}\times \Pb^N$ over $\Pb^N$ respectively, while the lower horizontal map $\pi$ is just a coordinate projection.   

Now let $U\subset \Ab_N(\bd')$ be the open subset of $\bf$ satisfying the conditions of~\eqref{eqn:def_of_fj}. The projection $\pi$ is an open map, so $V:=\pi(U)\subset \Ab_N(\bd)$ is open. We claim that for every $\bf\in V(K)$, $X_\bf(K)\subset X_\bf$ is dense. To see this, let
\[
    \pi^{-1}(\bf):=\{(\bg,\bf)\in\Ab_N(\bd')\}.
\]
Counting equations, we have $\dim X_{\bg,i}\ge N-(|\bm'|-|\bm|)$ (resp. $\dim X_{\bf,j}\ge N-|\bm|$) for each irreducible component $X_{\bg,i}$ of $X_{\bg}$ (resp. $X_{\bf,j}$ of $X_\bf$). By Lemma~\ref{lemma:nonempty}, $N\ge |\bm'|$, so we conclude that $X_{\bg,i}\cap X_{\bf,j}\neq \emptyset$ for each $i,j$; in particular every irreducible component of $X_{\bg}$ meets every irreducible component of $X_{\bf}$ and $X_{(\bg,\bf)}:=X_\bg\cap X_\bf\subset\Pb^N$ is non-empty. Note also that the restriction of the upper horizontal map above
\[
    \mcx_N(\bd')|_{X_{\bf}}\to X_{\bf}
\]
is a vector bundle. In particular, the irreducible components of $\mcx_N(\bd')|_{X_{\bf}}$ correspond bijectively to those of $X_\bf$, and each is a vector bundle over its image. 

Taken together, the above implies that the projection
\[
    \mcx_N(\bd')|_{X_{\bf}}\to \pi^{-1}(\bf)
\]
restricts to a surjection on every irreducible component of $\mcx_N(\bd')|_{X_{\bf}}$. As $(U\cap\pi^{-1}(\bf))(K)\subset \pi^{-1}(\bf)$ is dense (it's a non-empty open in an affine space), and as $X_{(\bg,\bf)}(K)\subset X_{(\bg,\bf)}$ is dense for all $(\bg,\bf)\in(U\cap \pi^{-1}(\bf))(K)$ (by assumption), we conclude by Lemma~\ref{lemma:fibration} that $\mcx_N(\bd')|_{X_\bf}(K)\subset \mcx_N(\bd')|_{X_\bf}$ is dense. The density of $X_{\bf}(K)\subset X_\bf$ follows.
\end{proof}
 
\subsection{Existence of solvable \texorpdfstring{$j$}{j}-planes}\label{subsection:existence}
We prove existence as a warm-up to density; the argument is the skeleton of the density argument of Section \ref{subsection:obliteration}. Note that existence holds for \emph{every} intersection of hypersurfaces in sufficiently many variables. 
 
\begin{definition}\label{definition:gj}
We define $g: \Nb \times \Nb^3 \to \Nb \cup \{\infty\}$ by
\[
\begin{aligned}
    \existence{j}{m_2,m_3,m_4}\coloneqq\min\{ N_0 \mid\ &\text{for every } N \geq N_0, \text{ for every } X \subseteq \Pb^N_{K} \\
    & \text{ of type } (0,m_2,m_3,m_4), \mcf(j,X)(K) \neq\emptyset \}.
\end{aligned}
\]
\end{definition}
 
Since $X$ ranges over all intersections defined over $K$, and a $k$-scheme $X$ has $\mcf(j,X)(K)=\mcf(j,X_K)(K)$, this also controls intersections defined over $k$ or over any intermediate field. Analogues for $\existence{j}{m_2,m_3,m_4}$ of the elementary properties of $\density{j}{m_2,m_3,m_4}$ are immediate:
 
\begin{lemma}\label{lemma:gj_elementary}
Let $j \in \Nb$ and $m_2, m_3, m_4 \in \Nb$.
\begin{enumerate}[label=(\alph*)]
    \item\label{i:g_base} We have $\existence{j}{0,0,0}=j$.
    \item\label{i:g_mono} If $m_i\leq m_i'$ for each $i \in \{2,3,4\}$, then $\existence{j}{m_2,m_3,m_4}\leq\existence{j}{m_2',m_3',m_4'}$.
    \item\label{i:g_linear} If $m_1\geq1$ and $N\geq\existence{j}{m_2,m_3,m_4}+m_1$, each $X\subseteq\Pb^N_K$ of type $(m_1,m_2,m_3,m_4)$ has $\mcf(j,X)(K)\neq\emptyset$.
\end{enumerate}
\end{lemma}
 
\begin{lemma}\label{lemma:planes_from_pts_existence}
Fix $j,m_2,m_3,m_4 \in \Nb$. Then
\[
    \existence{j}{m_2,m_3,m_4} \leq \existence{0}{m_2+j m_3+\tbinom{j+1}{2} m_4,\ m_3+j m_4,\ m_4} + j + j m_2 + \tbinom{j+1}{2} m_3+\tbinom{j+2}{3} m_4 .
\]
\end{lemma}
\begin{proof}
By Example~\ref{example:type_shift}, the three arguments of $\existence{0}{\cdot}$ on the right-hand side are the last three entries of $\bm^j$, and the remaining summand $jm_2+\tbinom{j+1}{2}m_3+\tbinom{j+2}{3}m_4$ is its first entry, the number of linear forms. Hence, by Lemma~\ref{lemma:gj_elementary}\ref{i:g_linear}, if $N$ is at least the right-hand side, then every intersection of type $\bm^j$ in $\Pb^{N-j}_K$ contains a $K$-point.

We induct on $j$, the case $j=0$ being trivial. Let $j\geq 1$, let $N$ be at least the right-hand side, and let $X\subseteq\Pb^N_K$ have type $\bm$. Since $\bm^{j-1} \leq \bm^j$ by Lemma~\ref{l:comb}\ref{i:norm}, monotonicity and the inductive hypothesis give $N\geq\existence{j-1}{m_2,m_3,m_4}$, and so $X$ contains a $K$-rational $(j-1)$-plane $\Lambda$. By Theorem~\ref{thm:type_of_fano}\ref{i:tf2}, $\mcf_\Lambda(j,X)^{\red}$ is isomorphic over $K$ to an intersection in $\Pb^{N-j}_K$ of type $\bm^j$. By assumption, this subscheme admits a $K$-point, i.e., there exists a $K$-rational $j$-plane $\Lambda'\supseteq\Lambda$ contained in $X$.
\end{proof}

In order to begin our induction, we need rational points on degree $\le 4$ subschemes of projective space; this is where $K$ being solvably closed enters.
 
\begin{lemma}\label{lemma:lowdegree}
Let $K$ be solvably closed, and let $W \subseteq \Pb^n_K$ be a closed subscheme which is either
\begin{enumerate}[label=(\alph*)]
\item the zero locus of a single form of degree at most $4$ (possibly the zero form), or
\item the zero locus of two quadratic forms, where $n \ge 2$.
\end{enumerate}
Then $W(K)$ is Zariski dense in $W$, and $W \neq \emptyset$ if $n \ge 1$.
\end{lemma}
\begin{proof}
Part (a) and non-emptiness are clear. For (b), suppose $W=\Vb(Q_1,Q_2)\subset\Pb^n_K$ with $n\ge2$. If $Q_1$ or $Q_2$ is zero, we are in case (a). Otherwise, we consider two cases. 

Firstly, if $Q_1$ and $Q_2$ share a common factor $h$, then we have
\[
  W=\Vb(h)\cup\Vb(Q_1/h,Q_2/h)
\]
set-theoretically. We see that $\Vb(h)$ is either a hyperplane or quadric, thus with dense $K$-points by (a), and hence $\Vb(Q_1/h,Q_2/h)$ is either linear or empty, respectively, and thus also with dense $K$-points. 

Secondly, $Q_1$ and $Q_2$ have no common factor, and every irreducible component $W'$ of $W$ is codimension $2$. The lines meeting $W$ form a proper closed subset of $\grass(1,n)_K$ of dimension at most \[\dim W+(n-1)=2n-3<\dim\grass(1,n).\] Hence, since $K$ is infinite, there is a $K$-rational line $\ell$ disjoint from $W$. Moreover, projecting $\pi_\ell:W'\to\Pb^{n-2}_K$ is finite and surjective (since $\dim W'=n-2$), so for a non-empty open $V \subseteq W'$ we know that $\pi_\ell(V)$ is dense and, in particular, contains a $K$-point $z$. The plane $\Lambda=\langle\ell,z\rangle$ is $K$-rational and $W\cap\Lambda=\pi_\ell^{-1}(z)$ is finite, so $Q_1|_\Lambda$ and $Q_2|_\Lambda$ have no common factor. Thus $W\cap\Lambda$ is a zero-dimensional scheme of length $4$ by B\'ezout, so each closed point $w$ has $[K(w):K]\le4$ and thus $K(w) = K$ by solvable closure. Since $z \in \pi_\ell(V)$, one of these points lies in $V$, hence every non-empty open contains a point in $W(K)$.
\end{proof}
 
\begin{lemma}\label{lemma:pts_from_planes_existence}
For any $m_2, m_3, m_4\in \Nb$, we have
\[ \existence{0}{m_2,m_3,m_4} \leq \min\left\{ \existence{1}{m_2-1,m_3,m_4}, \; \existence{2}{m_2-2,m_3,m_4}, \; \existence{1}{m_2,m_3-1,m_4}, \; \existence{1}{m_2,m_3,m_4-1} \right\}, \]
where terms involving negative arguments are omitted.
\end{lemma}
\begin{proof}
Fix one of the four cases, and let $X=X_\bf\subseteq\Pb^N_K$ have type $(0,m_2,m_3,m_4)$ with $N$ at least the corresponding term. Split $\bf=(\bf_Y,\bf_Z)$, where $\bf_Y$ consists of one quadric, two quadrics, one cubic, or one quartic respectively, and $\bf_Z$ consists of the remaining forms: of type $(0,m_2-1,m_3,m_4)$, $(0,m_2-2,m_3,m_4)$, $(0,m_2,m_3-1,m_4)$, or $(0,m_2,m_3,m_4-1)$ respectively. We set $Y\coloneqq \Vb(\bf_Y)$ and $Z\coloneqq \Vb(\bf_Z)$, so $X=Y\cap Z$, and let $j=1,2,1,1$ respectively, so that in each case $\bf_Y$ consists of $j$ forms. By the choice of $N$, there is a $K$-rational $j$-plane $\Lambda\subseteq Z$. Then $Y\cap\Lambda\subseteq\Lambda\cong\Pb^j_K$ is the zero locus of $j$ forms on $\Pb^j_K$ of the kind considered in Lemma~\ref{lemma:lowdegree}, hence has a $K$-point, which lies in $Y\cap Z=X$.
\end{proof}

\subsection{Obliteration lemmas}\label{subsection:obliteration}
The heart of our argument relies on two fundamental inequalities, which we prove here, and which essentially go back to Sylvester \cite{Sylvester1887}. We first consider the hypothesis of Corollary~\ref{cor:DM_locus}. Although the remainder of this section concerns types $(0,m_2,m_3,m_4)$, the lemma holds for intersections of arbitrary degrees.

\begin{lemma}\label{lemma:numerics}
Let $j \geq 1$ and let $\bm = (0, m_2, \dots, m_n) \neq \mathbf{0}$ be a type with no linear forms, with degree tuple $\bd = (d_1,\dots,d_r)$, and suppose $N \geq |\bm^j| + j$. Then $\delta^-(N,\bd,j) > 0$, except for the case of a single quadric and $N = 2j+1$, in which case $\delta^-(N,\bd,j) = 0$.
\end{lemma}
\begin{proof}
By Lemma~\ref{l:comb}\ref{i:norm}, $|\bm^j| = \sum_{i=1}^r \tbinom{d_i+j-1}{j}$, and $\tbinom{d+j}{j} = \tfrac{d+j}{d}\tbinom{d+j-1}{j}$. Using $N - j \geq |\bm^j|$, we compute
\[
\begin{aligned}
    \delta(N,\bd,j) & = (j+1)(N-j) - \sum_{i=1}^r \tbinom{d_i+j}{j} \\
    & \geq \sum_{i=1}^r \tbinom{d_i+j-1}{j}\left(j+1-\tfrac{d_i+j}{d_i}\right) \\
    & = \sum_{i=1}^r \tbinom{d_i+j-1}{j}\,\tfrac{j(d_i-1)}{d_i} > 0,
\end{aligned}
\]
since all $d_i \geq 2$ and $r \geq 1$. Next, since $\tbinom{d+j-1}{j} \geq j+1$ for $d \geq 2$, we have
\[
    N - 2j - r \;\geq\; |\bm^j| - j - r = \sum_{i=1}^r \left(\tbinom{d_i+j-1}{j} - 1\right) - j \;\geq\; (r-1)j \;\geq\; 0,
\]
with equality throughout if and only if $r = 1$, $d_1 = 2$, and $N = |\bm^j| + j = 2j+1$.
\end{proof}

\begin{lemma}[Linear subspaces via polar cones]\label{lemma:planes_from_pts}
Fix $j,m_2,m_3,m_4 \in \Nb$. Then
\[
    \density{j}{m_2,m_3,m_4} \leq \density{0}{m_2+j m_3+\tbinom{j+1}{2} m_4,\ m_3+j m_4,\ m_4} + j + j m_2 + \tbinom{j+1}{2} m_3+\tbinom{j+2}{3} m_4 .
\]
\end{lemma}
\begin{proof}
Write $R_j$ for the right-hand side of the inequality. As in the proof of Lemma~\ref{lemma:planes_from_pts_existence}, the arguments of $\density{0}{\cdot}$ in $R_j$ are the last three entries of $\bm^j$ and the remaining summand $R_j-j-\density{0}{\cdot}$ is its first entry, i.e. Lemma~\ref{lemma:omit_linear} applies to intersections of type $\bm^j$ in $\Pb^{N-j}_K$ when $N\geq R_j$. Moreover, Lemma~\ref{lemma:nonempty} and \eqref{eqn:floor} give
\begin{equation}\label{eqn:Rj_floor}
    R_j \geq |\bm^j|+j.
\end{equation}
Note that the inequality is vacuously true if $R_j=\infty$, so it suffices to assume $R_j<\infty$ for all $j$ for the purpose of proving the lemma. 
We will induct on $j$ to show the inequality above, the case $j=0$ being trivial.

Non-trivially, we suppose that $(m_2,m_3,m_4)\neq(0,0,0)$ and $R_j<\infty$. Let $j\geq1$ and fix $N\geq R_j$. By Definition~\ref{definition:fj} we must show that for general $\bf\in\param{N}{\bd}$, $\mcf(j,X_\bf)(K)$ is non-empty and dense in $\mcf(j,X_\bf)$ and that $p_{j,0}:\mcc^j(X_\bf)\to X_\bf$ is surjective.

By \eqref{eqn:Rj_floor}, Theorem~\ref{thm:type_of_fano} says that for \emph{every} intersection $X\subseteq\Pb^N$ of type $\bm$ (over any field), we have
\begin{equation}\label{eqn:pj_surjective}
\mcf_\Lambda(j,X)\neq\emptyset\ \text{ for every } \Lambda\in\mcf(j-1,X),\end{equation}
In particular, $p_j:\mcc^j(X)\to\mcc^{j-1}(X)$ and thus $p_{j,0}:\mcc^j(X)\to X$ are surjective maps (of varieties).

Similarly, Lemma~\ref{lemma:omit_linear} (with $j=0$, applied to the type $\bm^j$) gives a dense open $U_\Pi\subseteq\param{N-j}{\bd^j}$ such that
\begin{equation}\label{eqn:U_Pi}
    X_\bf(K)\ \text{is non-empty and dense in } X_\bf\subseteq\Pb^{N-j}_K\ \text{for all } \bf\in U_\Pi(K).
\end{equation}
Without loss of generality (i.e. replacing this by its $\mcg_{N-j}(\bd^j)$-orbit if necessary, as observed after Definition~\ref{definition:fj}), we can assume that $U_\Pi$ is $\mcg_{N-j}(\bd^j)$-stable.

Finally, since $\bm^{j-1}\leq\bm^j$ entrywise, Lemma~\ref{lemma:monotonicity} gives $R_{j-1}\leq R_j\leq N$ and so, by the inductive hypothesis $N\geq\density{j-1}{m_2,m_3,m_4}$, there is a dense open $U_{j-1}\subseteq\param{N}{\bd}$ such that 
\begin{equation}\label{eqn:Ujm1}
    \mcf(j-1,X_\bf)(K)\ \text{is dense in }\mcf(j-1,X_\bf)\ \text{for all }\bf\in U_{j-1}(K).
\end{equation}
The proof proceeds through constructing a dense open $U_j\subseteq\param{N}{\bd}$. Then, for $\bf\in U_j(K)$, we verify the hypotheses of Lemma~\ref{lemma:fibration} for $p_j:\mcc^j(X_\bf)\to\mcc^{j-1}(X_\bf)$ in order to lift dense supplies of points.

\medskip
We define three dense opens of $\param{N}{\bd}$, $U_{\mathrm{eq}}$, $U_{\mathrm{tr}}$, and $U_{\mathrm{irr}}$, and set 
\[
    U_j:=U_{j-1}\cap U_{\mathrm{eq}}\cap U_{\mathrm{tr}}\cap U_{\mathrm{irr}}.
\]
Firstly, we let $U_{\mathrm{eq}}$ be the open of Proposition~\ref{prop:equidim} for $j-1$: for $\bf\in U_{\mathrm{eq}}$, every irreducible component of $\mcf(j-1,X_\bf)$ has dimension $\delta_{j-1}=\delta(N,\bd,j-1)$. 

Next, consider the Zariski-locally trivial fibration
\[
    \mcf(j-1,\mcx)\to \grass(j-1,N)
\]
whose fiber over $\Lambda_0$ is $\param{N,\Lambda_0}{\bd}$, a linear subspace and hence irreducible. Consider
\[
    U_0\coloneqq\Pi_0^{-1}(U_\Pi)\subseteq\param{N,\Lambda_0}{\bd},
\]
which is a dense $\mcg_{\Lambda_0}$-stable open by Theorem~\ref{thm:linear_model}, parts~\ref{i:lm_linear} and~\ref{i:lm_equivariant}.  Since the $\mcg$-action on $\mcf(j-1,\mcx)$ induces a transitive action of $\mcg$ on $\grass(j-1,N)$, we obtain a dense open
\[
    \mcg\cdot U_0\subseteq\mcf(j-1,\mcx).
\]
Set $\mcz\coloneqq\mcf(j-1,\mcx)\setminus\mcg\cdot U_0$. This is a $\mcg$-invariant closed subset with 
\[
	\dim\mcz<\dim\mcf(j-1,\mcx)=\dim\param{N}{\bd}+\delta_{j-1}
\]
by Lemma~\ref{lemma:universal_smooth}. We apply Lemma~\ref{lemma:general-dim} to the $\mcg$-equivariant projection $\mcz\to\param{N}{\bd}$ with $\delta=\delta_{j-1}-1$, to obtain a dense open $U_{\mathrm{tr}}$ consisting of all $\bf\in \param{N}{\bd}$, with
\[
    \dim\left(\mcz\cap\mcf(j-1,X_\bf)\right) \leq \delta_{j-1}-1.
\]
Lastly, we define $U_{\mathrm{irr}}$ by the two cases of Lemma~\ref{lemma:numerics}: if $\delta^-(N,\bd,j)>0$ then we set $U_{\mathrm{irr}}\coloneqq\param{N}{\bd}\setminus\Delta_j$, with $\Delta_j$ as in Corollary~\ref{cor:DM_locus}; otherwise, we have $(m_2,m_3,m_4)=(1,0,0)$ and $N=2j+1$, and we take $U_{\mathrm{irr}}$ as the open of non-degenerate quadratic forms. 

We now claim that $U_j \coloneqq U_{j-1} \cap U_{\mathrm{eq}} \cap U_{\mathrm{tr}} \cap U_{\mathrm{irr}},$ has the property that, for all $\bf\in U_j(K)$, $\mcf(j,X_{\bf})(K)\subset\mcf(j,X_{\bf})$ is dense and $p_{j,0}: \mcc^j(X_{\bf})\to X_\bf$ is surjective. Because $U_j\subset U_{j-1}$, this will follow from the inductive hypothesis and Lemma~\ref{lemma:fibration}, once we verify that the hypotheses of Lemma~\ref{lemma:fibration} apply to the map $p_j: \mcc^j(X_\bf)\to \mcc^{j-1}(X_\bf)$.

\medskip
\noindent\textbf{Hypothesis \ref{item:dense-base}.} 
We take $V\coloneqq(\mcg\cdot U_0)\cap\mcf(j-1,X_\bf)$, an open subset of $\mcf(j-1,X)$. Since $\bf\in U_{\mathrm{eq}}\cap U_{\mathrm{tr}}$, every irreducible component of $\mcf(j-1,X)$ has dimension $\delta_{j-1}>\dim\left(\mcz\cap\mcf(j-1,X_\bf)\right)$, so no component is contained in $\mcz\cap\mcf(j-1,X_\bf)$ and $V$ is dense in $\mcf(j-1,X_\bf)$. Moreover, because $\bf\in U_{j-1}$, we know $\mcf(j-1,X_\bf)(K) \subseteq \mcf(j-1,X_\bf)$ is dense by \eqref{eqn:Ujm1} and so $V(K) \subseteq V$ is too. Recall the Zariski locally trivial bundle
\[
    q_{j-1}: \mcc^{j-1}(X_\bf)\to \mcf(j-1,X_\bf)
\]
of Proposition~\ref{prop:span_properties}\ref{i:qbundle}. Restricting to $V$, we obtain a Zariski locally trivial bundle with geometrically irreducible rational fiber, and so Corollary~\ref{corollary:density-bundle} shows that
\[
    \Sigma\coloneqq q_{j-1}^{-1}(V)(K)
\]
is dense in $q_{j-1}^{-1}(V)$, which is a dense open in $\mcc^{j-1}(X_\bf)$. Thus $\Sigma\subseteq \mcc^{j-1}(X_\bf)(K)$ is dense in $\mcc^{j-1}(X_\bf)$.

\medskip
\noindent\textbf{Hypothesis \ref{item:dense-fibers}.}
Let $b=(x_0,\dots,x_{j-1},\Lambda)\in \Sigma$, so that $\Lambda\in V(K)$, i.e., $(\bf,\Lambda)\in(\mcg\cdot U_0)(K)$. Since $\Lambda$ is $K$-rational, we may choose $g\in\GL_{N+1}(K)$ with $g\cdot\Lambda_0=\Lambda$ by extending a $K$-basis of the subspace underlying $\Lambda$ to a basis of $K^{N+1}$. Translating back, we have $g^{-1}\cdot\bf\in U_0(K)$, so we may set
\[
    \bf_\Lambda\coloneqq\Pi_0(g^{-1}\cdot\bf)\in U_\Pi(K)\qquad\text{and}\qquad X_{\bf_\Lambda}\subseteq\Pb^{N-j}_K ,
\]
the \emph{polar model} of $X_{g^{-1}\cdot\bf}$ at $\Lambda_0$. Since $X_{g^{-1}\cdot\bf}=g^{-1}(X_\bf)$ and $\Lambda_0=g^{-1}(\Lambda)$, the automorphism $g$ of $\Pb^N_K$ induces an isomorphism $\mcf_{\Lambda_0}(j,X_{g^{-1}\cdot\bf})\cong\mcf_\Lambda(j,X_\bf)$ of $K$-schemes, and Theorem~\ref{thm:linear_model}\ref{i:lm_iso} gives
\[
    \mcf_\Lambda(j,X)^{\red}\ \cong\ \mcf_{\Lambda_0}(j,X_{g^{-1}\cdot\bf})^{\red}\ \cong\ X_{\bf_\Lambda}.
\]
By \eqref{eqn:U_Pi}, $X_{\bf_\Lambda}(K)$ is dense in $X_{\bf_\Lambda}$, so $\mcf_\Lambda(j,X)^{\red}(K)$ is dense in $\mcf_\Lambda(j,X)$. By Proposition~\ref{prop:span_properties}\ref{i:pfiber}, $p_j^{-1}(b)^{\red}$ is a Zariski locally trivial $\Ab^j$-bundle over $\mcf_\Lambda(j,X)^{\red}$.  By \eqref{eqn:U_Pi} and Corollary~\ref{corollary:density-bundle}, we conclude that $p_j^{-1}(b)(K)\subset p_j^{-1}(b)$ is dense.

\medskip
\noindent\textbf{Hypothesis \ref{item:dominate-components}.}
We consider the two cases. Firstly, when $\delta^-(N,\bd,j)>0$, Corollary~\ref{cor:DM_locus} shows that $\mcf(j,X_\bf)$ is geometrically irreducible since $\bf\notin\Delta_j$. As $q_j: \mcc^j(X_\bf)\to\mcf(j,X_\bf)$ is a Zariski locally trivial bundle with geometrically irreducible fibers (Proposition~\ref{prop:span_properties}\ref{i:qbundle}), $\mcc^j(X_\bf)$ is irreducible. By the discussion around \eqref{eqn:pj_surjective}, $p_j$ is surjective, so $\mcc^{j-1}(X_\bf)$ is irreducible, and the unique component $\mcc^j(X_\bf)$ dominates it.

When $\delta^-(N,\bd,j)=0$, we have $(m_2,m_3,m_4)=(1,0,0)$ and $N=2j+1$, so $X=Q$ is a smooth quadric in $\Pb^{2j+1}_K$. Classically, $\mcf(j,Q)$ is the disjoint union of two geometrically irreducible components (the two rulings of $Q$), and every $(j-1)$-plane $\Lambda\subset Q$ is contained in exactly one $j$-plane from each ruling \cite[Theorem~22.14]{Harris1992}. Hence $\mcc^j(Q)$ has exactly two irreducible components, the preimages of the two rulings, and each of them maps onto a component of $\mcc^{j-1}(Q)$.  In either case, the hypothesis holds.

\medskip
Now, by Lemma~\ref{lemma:fibration}, $\mcc^j(X)(K)$ is dense in $\mcc^j(X)$. Since $q_j:\mcc^j(X)\to\mcf(j,X)$ is surjective, we see that $\mcf(j,X)(K)$ is dense in $\mcf(j,X)$; it is non-empty because $X\neq\emptyset$ (as $N\geq|\bm^j|\geq|\bm|$) and every point of $X$ lies on a $j$-plane by \eqref{eqn:pj_surjective}, which also gives the surjectivity of $p_{j,0}$. As $U_j$ is a dense open and $\bf\in U_j(K)$ was arbitrary, a general intersection of type $\bm$ in $\Pb^N_K$ satisfies the two properties of \eqref{eqn:def_of_fj}. Since $N\geq R_j$ was arbitrary, Definition~\ref{definition:fj} gives $\density{j}{m_2,m_3,m_4}\leq R_j$, which completes the induction.
\end{proof}

\begin{lemma}[Obliteration via linear subspaces]\label{lemma:pts_from_planes}
For any $m_2, m_3, m_4\in \Nb$, we have
\[ \density{0}{m_2,m_3,m_4} \leq \min\left\{ \density{1}{m_2-1,m_3,m_4}, \; \density{2}{m_2-2,m_3,m_4}, \; \density{1}{m_2,m_3-1,m_4}, \; \density{1}{m_2,m_3,m_4-1} \right\}, \]
where terms involving negative arguments are omitted.
\end{lemma}
\begin{proof}
Fix one of the four cases, let $j=1,2,1,1$ respectively, and let $(m_2',m_3',m_4')$ be the arguments of the corresponding term, so that the term is $\density{j}{m_2',m_3',m_4'}$. As in the proof of Lemma~\ref{lemma:pts_from_planes_existence}, after reordering we write the degree tuple of $\bm$ as $\bd=\bd_Y\oplus\bd_Z$, where $\bd_Y$ is $(2)$, $(2,2)$, $(3)$ or $(4)$ respectively---so that $\bd_Y$ has exactly $j$ entries---and $\bd_Z$ is the degree tuple of $\bm_Z\coloneqq(0, m_2',m_3',m_4')$. Accordingly,
\[
    \param{N}{\bd}=\param{N}{\bd_Y}\times\param{N}{\bd_Z},
\]
and for $\bf=(\bf_Y,\bf_Z)$ we write $Y\coloneqq\Vb(\bf_Y)$ and $Z\coloneqq\Vb(\bf_Z)$, so that $X_\bf=Y\cap Z$. Note that this decomposition is \emph{not} preserved by $\mcr_N(\bd)$, which may add multiples of the forms of $\bf_Z$ to those of $\bf_Y$.
 
Suppose $\density{j}{m_2',m_3',m_4'}<\infty$ and fix $N\geq\density{j}{m_2',m_3',m_4'}$. We must show that a general intersection of type $\bm = (0, m_2, m_3, m_4)$ in $\Pb^N_K$ has non-empty dense $K$-points. 
 
\medskip
For $\bf=(\bf_Y,\bf_Z)\in\param{N}{\bd}(K)$, consider $\mcf(j,Z)$ and $\varphi: S_j|_{\mcf(j,Z)}\to \mcf(j,Z)$ be the restriction of the tautological bundle regarded as the incidence variety $\{(y,\Lambda)\mid y\in\Lambda\in \mcf(j,Z)\}\subseteq\Pb^N_K\times \mcf(j,Z)$, which is a Zariski locally trivial $\Pb^j$-bundle over $\mcf(j,Z)$. We write $E\coloneqq\{ (y,\Lambda) \in S_j|_{\mcf(j,Z)} \mid y \in Y \}$, $\pi\coloneqq\varphi|_E: E\to \mcf(j,Z)$, and 
\begin{align*}
    p: E &\to Y\cap Z\\
    (y,\Lambda) &\mapsto y,
\end{align*}
so that we have the correspondence
\[ \begin{tikzcd}
& E \arrow[swap]{dl}{p} \arrow{dr}{\pi} \\
Y\cap Z && \mcf(j,Z).
\end{tikzcd} \]
Note that $E$ is the common zero locus in $S_j|_{\mcf(j,Z)}$ of the $j$ sections of $\mco_{S_j}(d)$ for $d\in\bd_Y$, obtained by restricting the forms of $\bf_Y$ to the tautological bundle; that $\pi^{-1}(\Lambda)=Y\cap\Lambda\subseteq\Lambda\cong\Pb^j$; and that $\pi(p^{-1}(y))\cong\mcf_y(j,Z)$.
 
Now, let $U_Z\subseteq\param{N}{\bd_Z}$ be a dense open on which the two properties of \eqref{eqn:def_of_fj} hold for $j$, which exists since $N\geq\density{j}{m_2',m_3',m_4'}$, and let $U_{\mathrm{eq}}\subseteq\param{N}{\bd_Z}$ be the open of Proposition~\ref{prop:equidim}, on which every irreducible component of $\mcf(j,Z)$ has dimension $\delta_Z\coloneqq\delta(N,\bd_Z,j)$.
 
To control the fibers of $\pi$, for $s\geq1$, let
\[
    D_{\geq s}(\Lambda)\ \subseteq\ \prod_{d\in\bd_Y}H^0\big(\Lambda,\mco(d)\big)
\]
be the locus of tuples of forms on $\Lambda$ whose common zero locus has dimension at least $s$. It is closed, by upper semicontinuity of fiber dimension for the proper incidence correspondence $\{(\bh,x)\mid \bh(x)=0\}$, and it is $\GL_{\Lambda}$-stable. Its universal version is
\[
    \mathcal D_{\geq s}\ \coloneqq\ \left\{(\bf_Y,\bf_Z,\Lambda)\in\param{N}{\bd_Y}\times\mcf(j,\mcx_Z) \mid \bf_Y|_\Lambda\in D_{\geq s}(\Lambda)\right\},
\]
where $\mcx_Z\coloneqq\mcx_N(\bd_Z)$ is the universal intersection of degree $\bd_Z$ and $\bf_Y|_\Lambda$ denotes the restriction of the forms of $\bf_Y$ to $\Lambda$. Restriction of forms to $\Lambda$ is a surjective linear map, varying algebraically with $\Lambda$ and trivializable Zariski locally on $\mcf(j,\mcx_Z)$, so $\mathcal D_{\geq s}$ is a closed subset of $\param{N}{\bd_Y}\times\mcf(j,\mcx_Z)$, of codimension $\codim D_{\geq s}$, the codimension of $D_{\ge s}(\Lambda)$ in $\prod_{d\in \bd_Y} H^0(\Lambda,\mco(d))$, which is independent of $\Lambda$. By Lemma~\ref{lemma:universal_smooth},
\[
    \dim\mathcal D_{\geq s}\ \leq\ \dim\param{N}{\bd}+\delta_Z-\codim D_{\geq s} .
\]
The fiber of $\mathcal D_{\geq s}$ over $\bf\in\param{N}{\bd}$, under the projection forgetting $\Lambda$, is
\[
    B_{\geq s,\bf}=\{\Lambda\in\mcf(j,Z)\mid\dim(Y\cap\Lambda)\geq s\}=\{\Lambda\in \mcf(j,Z)\mid\dim\pi^{-1}(\Lambda)\geq s\}.
\]
By Lemma~\ref{lemma:general-dim} there is a dense open $U_s\subseteq\param{N}{\bd}$ such that $\dim B_{\geq s,\bf}\leq\delta_Z-\codim D_{\geq s}$ for all $\bf\in U_s$. We consider the dense open subset of $\param{N}{\bd}$ given by
\[
    U\ \coloneqq\ \big(\param{N}{\bd_Y}\times(U_Z\cap U_{\mathrm{eq}})\big)\ \cap\ \bigcap_{s=1}^{j}U_s.
\]
Fix $\bf\in U(K)$; we verify the hypotheses of Lemma~\ref{lemma:fibration} for $\pi: E\to \mcf(j,Z)$.

\medskip
\noindent\textbf{Hypothesis \ref{item:dense-base}.} Take $\Sigma\coloneqq \mcf(j,Z)(K)$, which is dense in $\mcf(j,Z)$ since $\bf_Z\in U_Z(K)$.
 
\medskip
\noindent\textbf{Hypothesis \ref{item:dense-fibers}.} For $\Lambda\in\Sigma$, the fiber $\pi^{-1}(\Lambda)=Y\cap\Lambda\subseteq\Lambda\cong\Pb^j_K$ is the zero locus of one form of degree $2$, $3$ or $4$ on $\Pb^1_K$, or of two quadratic forms on $\Pb^2_K$ (by our assumptions on $j$ above). By Lemma~\ref{lemma:lowdegree}, its $K$-points are dense, and it is non-empty. 
 
\smallskip
\noindent\textbf{Hypothesis \ref{item:dominate-components}.} We apply Lemma~\ref{lemma:dominance} to $E\subseteq S_j|_{\mcf(j,Z)}\xrightarrow{\varphi}\mcf(j,Z)$ with $a=c=j$. For every irreducible component $B'$ of $\mcf(j,Z)$, $\varphi^{-1}(B')$ is a $\Pb^j$-bundle over $B'$, hence irreducible of dimension $\dim B'+j$, and $E$ is cut out of $S_j|_{\mcf(j,Z)}$ by $c=j$ sections of line bundles. Since $\bf_Z\in U_{\mathrm{eq}}$, every component $B'$ has $\dim B'=\delta_Z$, so for $s>a-c=0$:
\[
    \dim(B_{\geq s,\bf}\cap B')+s \leq \dim B_{\geq s,\bf}+s \leq \delta_Z-\codim D_{\geq s}+s = \dim B'+(a-c)-(\codim D_{\geq s}-s),
\]
and the hypothesis \eqref{eq:domcount} of Lemma~\ref{lemma:dominance} holds provided $\codim D_{\geq s}>s$ for every $1\leq s\leq j$. We check this in each case:
\begin{itemize}
    \item $j=1$, $\bd_Y=(d)$ with $d\in\{2,3,4\}$, $s=1$: a form on $\Pb^1$ vanishes on a positive-dimensional set only if it is zero, so $D_{\geq1}=\{0\}$ has codimension $d+1\geq3>1$.
    \item $j=2$, $\bd_Y=(2,2)$, $s=1$: two ternary quadratic forms $(Q_1,Q_2)$ have a common zero locus of dimension $\geq1$ if and only if they have a common factor over $\bar K$. The pairs with a common factor are: those with $Q_1=0$ or $Q_2=0$ (codimension $6$); non-zero proportional pairs (dimension $6+1=7$, codimension $5$); and pairs $(\ell\ell_1,\ell\ell_2)$ with a common linear factor $\ell$ (dimension at most $2+3+3=8$, codimension at least $4$). Hence $\codim D_{\geq1}=4>1$.
    \item $j=2$, $\bd_Y=(2,2)$, $s=2$: $D_{\geq2}=\{(0,0)\}$ has codimension $12>2$.
\end{itemize}
So every irreducible component of $E$ dominates an irreducible component of $\mcf(j,Z)$.
 
\medskip
By Lemma~\ref{lemma:fibration}, $E(K)$ is dense in $E$. The map $p: E\to Y\cap Z$ is surjective: a point $y\in Y\cap Z$, over any extension of $K$, lies on a $j$-plane $\Lambda\subseteq Z$ by the second property of \eqref{eqn:def_of_fj} for $\bf_Z\in U_Z(K)$, and then $(y,\Lambda)\in E$. Hence $(Y\cap Z)(K)=X_\bf(K)$ is dense in $X_\bf$, and it is non-empty because $E\neq\emptyset$: $\mcf(j,Z)\neq\emptyset$ by the first property of \eqref{eqn:def_of_fj}, and every fiber of $\pi$ is non-empty, as shown in the verification of Hypothesis~\ref{item:dense-fibers}. We conclude that a general intersection of type $\bm$ in $\Pb^N_K$ has non-empty dense $K$-points; and as $N\geq\density{j}{m_2',m_3',m_4'}$ was arbitrary, $\density{0}{m_2,m_3,m_4}\leq\density{j}{m_2',m_3',m_4'}$.
\end{proof}

\subsection{Iterative formulas for \texorpdfstring{$\density{j}{m_2,m_3,m_4}$}{fj(m2,m3,m4)}}\label{subsection:iterative_formulas}

With our workhorse lemmas established, all that remains is realizing the consequences of our inequalities.
For brevity, we will write
\[ \density{j}{m_2} \coloneqq \density{j}{m_2,0,0} \qquad \text{ and } \qquad \density{j}{m_2,m_3} \coloneqq \density{j}{m_2,m_3,0}, \]
and similarly for $\existence{j}{\cdot}$. Note that the arithmetic here applies verbatim to $\existence{j}{\cdot}$ in place of $\density{j}{\cdot}$.
The algorithm is straightforward: we repeatedly apply Lemmas \ref{lemma:pts_from_planes} and \ref{lemma:planes_from_pts} in order to ``obliterate'' quartics, then cubics, and then quadrics systematically until the base case of Lemma~\ref{lemma:base_case} is reached, namely
\[ \density{j}{0} = j \qquad \text{and} \qquad \existence{j}{0} = j. \]
Note that applying Lemma \ref{lemma:pts_from_planes} once, together with Lemma~\ref{lemma:nonempty}, recovers the basic facts
\[ \density{0}{1}=\density{0}{0,1}=\density{0}{0,0,1}=1 \qquad \text{ and } \qquad \density{0}{2} = 2, \]
which reflect the classical equations in radicals for quadric, cubic, and quartic polynomials. If we count these among our base cases, in practice, $\density{0}{1}$ and $\density{0}{2}$ are the terminal cases for most computations. All values that occur are finite by induction.

\begin{example}\label{example:point_on_quadric_cubic}
    Our results specialize to the classical fact (dating at least to \cite{Bring1786}) that a general intersection of a quadric $Q$ and a cubic $C$ in $\Pb^3$---indeed, as we will see, every smooth one---has a dense collection of solvable points. Indeed, the classical theory of quadrics shows that (a) $Q(k^{\Sol})\subset Q$ and $\mcf(1,Q)(k^{\Sol})\subset \mcf(1,Q)$ are dense so long as $\characteristic(k)\neq 2$ and $Q$ has a smooth point, and (b) every point lies on a line, i.e., $\density{1}{1}=3$. Picking a line $\ell\in\mcf(1,Q)(k^{\Sol})$, the solvability of the cubic in one variable shows that $\ell\cap C\subset (Q\cap C)(k^{\Sol})$. As in the proof of Lemma~\ref{lemma:pts_from_planes}, this process generates a dense collection of solvable points, provided hypothesis \ref{item:dominate-components} of Lemma~\ref{lemma:fibration} holds for the correspondence $E\to\mcf(1,Q)$. This is true when $Q\cap C$ contains no line, since then every fiber $\ell\cap C$ is finite and Lemma~\ref{lemma:dominance} applies with no excess locus (Remark~\ref{remark:dominance_krull}). This aligns with the computation:
    \[
        \density{0}{1,1}\le \density{1}{1}\le \density{0}{1}+2=3.
    \]
\end{example}


We argue the general cases by induction.

\begin{prop}\label{prop:points_on_quadrics} For all $m_2 \in \Nb$,
\[ \density{0}{m_2} \leq \floor{\tfrac{m_2+1}{2}}^2 + \floor{\tfrac{m_2}{2}}^2. \]
\end{prop}
\begin{proof}
The base case $\density{0}{0} = 0$ is trivial and $\density{0}{1} = 1$ is the quadratic formula. For any $m_2 \geq 2$, we have
\[ \begin{split}
    \density{0}{m_2} \leq \density{2}{m_2-2} & \leq \density{0}{m_2-2}+2(m_2-1) \\
    & \leq \floor{\tfrac{m_2-1}{2}}^2 + \floor{\tfrac{m_2-2}{2}}^2 + 2m_2 - 2 \\
    & = \left\{ \begin{array}{ll}
        (a-1)^2 + (a-1)^2 + 4a - 2, & \text{if } m_2 = 2a \\
        a^2 + (a-1)^2 + 4a, & \text{if } m_2 = 2a+1 \\
\end{array} \right. \\
    & = \left\{ \begin{array}{ll}
        a^2 + a^2, & \text{if } m_2 = 2a \\
        (a+1)^2 + a^2, & \text{if } m_2 = 2a+1 \\
\end{array} \right. \\
    & = \floor{\tfrac{m_2+1}{2}}^2 + \floor{\tfrac{m_2}{2}}^2.
\end{split} \]
\end{proof}

\begin{remark}
Obliterating quadrics in pairs is more efficient than doing so one at a time. If $m_2 \geq 2$, then
\[
    \density{0}{m_2} \leq \density{1}{m_2-1} \leq \density{0}{m_2-1} + m_2 \leq \density{1}{m_2-2} + m_2 \leq \density{0}{m_2-2} + 2m_2-1
\]
is a worse bound than $\density{0}{m_2} \leq \density{2}{m_2-2} \leq \density{0}{m_2-2} + 2m_2-2$.
\end{remark}

Together with Lemma \ref{lemma:planes_from_pts}, this computation generalizes a classical result of obtaining $j$-planes on quadric hypersurfaces by diagonalizing the quadric form:
\[
    \density{j}{1} \leq \density{0}{1} + 2j = 2j+1.
\]

\begin{prop}\label{prop:points_on_quadrics_anm_cubics} For all $m_2, m_3 \in \Nb$,
\[ \density{0}{m_2,m_3} \leq \density{0}{m_2+\tbinom{m_3}{2}}+\tfrac{1}{3}m_3(3m_2+m_3^2+2). \]
\end{prop}
\begin{proof}
We induct on $m_3$, where $m_3=0$ is tautological. For the inductive step, fix $m_2 \in \Nb$ and compute:
\begin{align*}
\density{0}{m_2,m_3} &\leq \density{1}{m_2,m_3-1}\\
& \leq \density{0}{m_2+m_3-1,m_3-1}+m_2+m_3 \\
& \leq \density{0}{m_2+m_3-1+\tbinom{m_3-1}{2}} + \tfrac{1}{3} (m_3-1)(3(m_2+m_3-1)+(m_3-1)^2+2) + m_2 + m_3 \\
& \leq \density{0}{m_2+\tbinom{m_3}{2}}+\tfrac{1}{3}m_3(3m_2+m_3^2+2).
\end{align*}
\end{proof}

\begin{prop}\label{prop:points_on_quadrics_and_cubics_and_quartics} For all $m_2, m_3, m_4 \in \Nb$,
\[ \density{0}{m_2,m_3,m_4} \leq \density{0}{m_2+m_3m_4+2\tbinom{m_4+1}{3},m_3+\tbinom{m_4}{2}} + m_2m_4 + m_3 \tbinom{m_4+1}{2} + \tfrac{1}{8}m_4(m_4+3)(m_4^2-m_4+2). \]
\end{prop}
\begin{proof}
Mirroring the proof of Proposition \ref{prop:points_on_quadrics_anm_cubics} we induct on $m_4$, where the base case is tautological. Applying the inductive hypothesis is straightforward:
\[
\begin{split}
\density{0}{m_2,m_3,m_4} & \leq \density{1}{m_2,m_3,m_4-1} \leq \density{0}{m_2+m_3+m_4-1,m_3+m_4-1,m_4-1}+m_2+m_3+m_4 \\
& \leq \density{0}{(m_2+m_3+m_4-1)+(m_3+m_4-1)(m_4-1)+2\tbinom{m_4}{3},(m_3+m_4-1)+\tbinom{m_4-1}{2}} \\
& \qquad + (m_2+m_3+m_4-1)(m_4-1) + (m_3+m_4-1)\tbinom{m_4}{2} \\
& \qquad \qquad + \tfrac{1}{8}(m_4-1)(m_4+2)((m_4-1)^2-(m_4-1)+2) + (m_2+m_3+m_4) \\
& \leq \density{0}{m_2+m_3m_4+2\tbinom{m_4+1}{3},m_3+\tbinom{m_4}{2}} + m_2m_4 + m_3 \tbinom{m_4+1}{2} + \tfrac{1}{8}m_4(m_4+3)(m_4^2-m_4+2).
\end{split}
\]
\end{proof}

We can now conclude by proving Theorems~\ref{theorem:solvable_points} and~\ref{theorem:solvable_j-planes}.

\begin{proof}[Proof of Theorem~\ref{theorem:solvable_j-planes}]
Combining Proposition \ref{prop:points_on_quadrics_and_cubics_and_quartics} with a single obliteration (Lemma \ref{lemma:planes_from_pts}) gives a bound for $f_j(m_2,m_3,m_4)$. Replacing the bound on $f_0(m_2)$ with the slightly worse $(\tfrac{m_2+1}{2})^2 + (\tfrac{m_2}{2})^2$ and unpacking Proposition \ref{prop:points_on_quadrics_and_cubics_and_quartics} in terms of Proposition \ref{prop:points_on_quadrics_anm_cubics} gives the polynomial $q(j,m_2,m_3,m_4)$ with the desired properties, written out in Appendix \ref{section:select_bounds}. All told, we have proven that
\[ f_j(m_2,m_3,m_4) \leq f_0\!\left(m_2+jm_3+\tbinom{j+1}{2}m_4,m_3+jm_4,m_4\right) + j + jm_2 + \tbinom{j+1}{2}m_3+\tbinom{j+2}{3}m_4 \leq q(j,m_2,m_3,m_4). \]
\end{proof}

\begin{proof}[Proof of Theorem~\ref{theorem:solvable_points}]
    This follows from Theorem~\ref{theorem:solvable_j-planes} by setting $j=0$ and $p(m_2,m_3,m_4)=q(0,m_2,m_3,m_4)$.
\end{proof}

\appendix

\section{Select bounds for \texorpdfstring{$f_j(m_2,m_3,m_4)$}{fj(m2,m3,m4)}}\label{section:select_bounds}

The polynomial $q(j,m_2,m_3,m_4)$ obtained through unpacking the Propositions of Section~\ref{subsection:iterative_formulas} is:

\begin{align*}
    q(j,m_2,m_3,m_4) & = \frac{1}{128} m_4^8 + \frac{1}{16} j m_4^7 + \frac{1}{96} m_4^7 + \frac{3}{16} j^2 m_4^6 + \frac{5}{48} j m_4^6 + \frac{1}{16} m_3 m_4^6 + \frac{17}{576} m_4^6 + \frac{1}{4} j^3 m_4^5 \\
    & \qquad + \frac{17}{48} j^2 m_4^5 + \frac{11}{48} j m_4^5 + \frac{5}{48} m_3 m_4^5 + \frac{3}{8} j m_3 m_4^5 + \frac{1}{48} m_4^5 + \frac{1}{8} j^4 m_4^4 + \frac{1}{2} j^3 m_4^4 + \frac{29}{48} j^2 m_4^4 \\
    & \qquad + \frac{3}{16} m_3^2 m_4^4 + \frac{11}{48} j m_4^4 + \frac{1}{8} m_2 m_4^4 + \frac{3}{4} j^2 m_3 m_4^4 + \frac{1}{6} m_3 m_4^4 + \frac{17}{24} j m_3 m_4^4 + \frac{169}{1152} m_4^4 \\
    & \qquad + \frac{1}{4} j^4 m_4^3 + \frac{7}{12} j^3 m_4^3 + \frac{31}{48} j^2 m_4^3 + \frac{7}{24} m_3^2 m_4^3 + \frac{3}{4} j m_3^2 m_4^3 + \frac{17}{24} j m_4^3 + \frac{1}{12} m_2 m_4^3 \\
    & \qquad + \frac{1}{2} j m_2 m_4^3+\frac{1}{2} j^3 m_3 m_4^3+\frac{3}{2} j^2 m_3 m_4^3+\frac{3}{16} m_3 m_4^3 + \frac{23}{24} j m_3 m_4^3 + \frac{3}{32} m_4^3 + \frac{1}{8} j^4 m_4^2 \\
    & \qquad + \frac{1}{2} j^3 m_4^2 + \frac{1}{4} m_3^3 m_4^2 + \frac{23}{24} j^2 m_4^2 + \frac{3}{4} j^2 m_3^2 m_4^2+\frac{5}{16} m_3^2 m_4^2+\frac{5}{4} j m_3^2 m_4^2+\frac{2}{3} j m_4^2 \\
    & \qquad + \frac{1}{2} j^2 m_2 m_4^2+\frac{3}{8} m_2 m_4^2+\frac{1}{2} j m_2 m_4^2+j^3 m_3 m_4^2+\frac{3}{2} j^2 m_3 m_4^2+\frac{25}{48} m_3 m_4^2+\frac{25}{24} j m_3 m_4^2 \\ 
    & \qquad +\frac{1}{2} m_2 m_3 m_4^2 + \frac{91}{288} m_4^2 + \frac{1}{6} j^3 m_4 + \frac{1}{4} m_3^3 m_4 + \frac{1}{2} j m_3^3 m_4+\frac{3}{4} j^2 m_4 +\frac{5}{4} j^2 m_3^2 m_4 \\
    & \qquad + \frac{5}{24} m_3^2 m_4+j m_3^2 m_4+j m_4+\frac{1}{2} j^2 m_2 m_4+j m_2 m_4+\frac{5}{12} m_2 m_4 +\frac{1}{2} j^3 m_3 m_4 \\
    & \qquad + \frac{5}{4} j^2 m_3 m_4+\frac{17}{12} j m_3 m_4+\frac{1}{2} m_2 m_3 m_4+j m_2 m_3 m_4+\frac{11}{24} m_3 m_4 +\frac{3}{8}  m_4 +\frac{1}{8} m_3^4 \\
    & \qquad + \frac{1}{2} j m_3^3 +\frac{1}{12} m_3^3 +\frac{1}{2} m_2^2 +\frac{1}{2} j^2 m_3^2+\frac{1}{2} j m_3^2 + \frac{1}{2} m_2 m_3^2+\frac{3}{8} m_3^2 + j + j m_2+\frac{1}{2} m_2 \\
    & \qquad + \frac{1}{2} j^2 m_3 + j m_3 + j m_2 m_3+\frac{m_2 m_3}{2}+\frac{5 m_3}{12}+\frac{1}{4},
\end{align*}
where the floor functions in Proposition \ref{prop:points_on_quadrics} have been dropped (a modest loss). In order to recover the bound in Theorem \ref{theorem:solvable_points}, we set $p(m_2,m_3,m_4) = q(0,m_2,m_3,m_4)$. We note that $p(m_2,m_3,0)$ has total degree $4$ and $p(m_2,0,0)$ is quadratic. A more natural printing of the terms of $p$ and $q$ is unknown to the authors. 

We also include tables of our best bounds on $\density{j}{m_2,m_3,m_4}$ and $\existence{j}{m_2,m_3,m_4}$ (which coincide) for the cases mentioned in Corollary \ref{corollary:points_on_pure_intersections}, and include the degree of the associated intersection $X$ for comparison. We emphasize that, in the first two tables, the bound on $N$ to guarantee dense solvable points is less than $\deg X$, let alone $2^{2^{\deg X}}$ as in \cite{Wooley2016}. Using the fact that exponential growth dominates polynomial growth for large $m_2, m_3,$ and $m_4$, we have
\[ 
    \density{0}{m_2,m_3,m_4} \leq p(m_2,m_3,m_4) < 2^{m_2} 3^{m_3} 4^{m_4} = \deg X 
\]
away from finitely many values of $m_2, m_3,$ and $m_4$. 
These estimates should also be compared with \cite[Theorem~2]{Wooley1998} which gives $N \geq 2m_2^2$, $N \geq 4m_3^4$, and $N \geq 16m_4^8$ to guarantee a solution to a system of $m_2$ quadrics, of $m_3$ cubics, and of $m_4$ quartics, respectively, or for a sharper comparison, with the formula obtained by solving the recursion implicit in \cite[Lemma 2.1]{Wooley1998} for $\varphi=1$.

\vspace{5mm} 

\begin{center}
\begin{tabular}{r|c|ccccccccc}
   & & & & & & $j$ \\
   $m_2$ & degree & 0 & 1 & 2 & 3 & 4 & 5 & 6 & 7 & 8 \\
   \hline
   1 & 2 & 1 & 3 & 5 & 7 & 9 & 11 & 13 & 15 & 17 \\
   2 & 4 & 2 & 5 & 8 & 11 & 14 & 17 & 20 & 23 & 26 \\
   3 & 8 & 5 & 9 & 13 & 17 & 21 & 25 & 29 & 33 & 37 \\
   4 & 16 & 8 & 13 & 18 & 23 & 28 & 33 & 38 & 43 & 48 \\
   5 & 32 & 13 & 19 & 25 & 31 & 37 & 43 & 49 & 55 & 61 \\
   6 & 64 & 18 & 25 & 32 & 39 & 46 & 53 & 60 & 67 & 74 \\
   7 & 128 & 25 & 33 & 41 & 49 & 57 & 65 & 73 & 81 & 89 \\
   8 & 256 & 32 & 41 & 50 & 59 & 68 & 77 & 86 & 95 & 104 \\
\end{tabular}
\captionof{table}{Bounds on \texorpdfstring{$f_j(m_2)$}{fj(m2)}} 
\end{center}

\vspace{5mm}

\begin{center}
\begin{tabular}{r|c|ccccccccc}
   & & & & & & $j$ \\
   $m_3$ & degree & 0 & 1 & 2 & 3 & 4 & 5 & 6 & 7 & 8 \\
   \hline
   1 & 3 & 1 & 5 & 10 & 18 & 27 & 39 & 52 & 68 & 85 \\
   2 & 9 & 5 & 16 & 33 & 56 & 85 & 120 & 161 & 208 & 261 \\
   3 & 27 & 16 & 42 & 81 & 131 & 194 & 268 & 355 & 453 & 564 \\
   4 & 81 & 42 & 95 & 168 & 261 & 374 & 507 & 660 & 833 & 1026 \\
   5 & 243 & 95 & 189 & 312 & 466 & 649 & 863 & 1106 & 1380 & 1683 \\
   6 & 729 & 189 & 340 & 533 & 768 & 1045 & 1364 & 1725 & 2128 & 2573 \\
   7 & 2187 & 340 & 568 & 853 & 1193 & 1590 & 2042 & 2551 & 3115 & 3736 \\
   8 & 6561 & 568 & 897 & 1298 & 1771 & 2316 & 2933 & 3622 & 4383 & 5216 \\
\end{tabular}
\captionof{table}{Bounds on \texorpdfstring{$f_j(0,m_3)$}{fj(0,m3)}} 
\end{center}

\vspace{5mm}

\begin{center}
\begin{tabular}{r|c|ccccccccc}
   & & & & & & $j$ \\
   $m_4$ & degree & 0 & 1 & 2 & 3 & 4 & 5 & 6 & 7 & 8 \\
   \hline
   1 & 4 & 1 & 10 & 44 & 133 & 319 & 656 & 1210 & 2059 & 3293 \\
   2 & 16 & 10 & 114 & 502 & 1476 & 3442 & 6918 & 12526 & 21000 & 33178 \\
   3 & 64 & 114 & 858 & 3180 & 8460 & 18510 & 35574 & 62328 & 101880 & 157770 \\
   4 & 256 & 858 & 4463 & 14028 & 33933 & 69758 & 128283 & 217488 & 346553 & 525858 \\
   5 & 1024 & 4463 & 17714 & 48650 & 108401 & 210797 & 372368 & 612344 & 952655 & 1417931 \\
   6 & 4096 & 17714 & 57680 & 142062 & 295218 & 546802 & 931756 & 1490318 & 2268014 & 3315666 \\
   7 & 16384 & 57680 & 161736 & 364492 & 713828 & 1267032 & 2090800 & 3261236 & 4863852 & 6993568 \\
   8 & 65536 & 161736 & 403665 & 845322 & 1573467 & 2690412 & 4314021 & 6577710 & 9630447 & 13636752 \\
\end{tabular}
\captionof{table}{Bounds on \texorpdfstring{$f_j(0,0,m_4)$}{fj(0,0,m4)}} 
\end{center}

\bibliography{bibliography}
\bibliographystyle{amsalpha}
\end{document}